\documentclass[a4paper,12pt]{article}
\usepackage{amssymb}
\usepackage{amsmath}
\usepackage{amsthm}
\usepackage{graphicx}
\numberwithin{equation}{section}

\newcommand{\ep}{\varepsilon}

\newcommand{\va}{\varphi}
\newcommand{\ppp}{\partial}
\newcommand{\dalpha}{\ppp_x^{\alpha}}

\newcommand{\weight}{e^{2s\va}}
\newcommand{\sumij}{\sum_{i,j=1}^n}
\newcommand{\sumalpha}{\sum_{\vert \alpha\vert \le 2}}
\newcommand{\sumalalpha}{\sum_{\vert \alpha\vert = 2}}

\newcommand{\R}{\mathbb{R}}
 
\newcommand{\N}{\mathbb{N}} 
\newcommand{\www}{\widetilde}
\newcommand{\OOO}{\Omega}
\newcommand{\ooo}{\overline}

\allowdisplaybreaks

\title{Carleman estimate and application to an inverse source problem 
for a viscoelasticity model in anisotropic case}

\author{Paola Loreti 
\thanks{
Dipartimento di Scienze di Base e Applicate per l'Ingegneria, Sapienza 
Universit\`a di Roma, Via Antonio Scarpa 16, 00161 Roma, Italy 
e-mail: {\tt paola.loreti@sbai.uniroma1.it}},\
Daniela Sforza 
\thanks{
Dipartimento di Scienze di Base e Applicate per l'Ingegneria, Sapienza 
Universit\`a di Roma, Via Antonio Scarpa 16, 00161 Roma, Italy 
e-mail: {\tt daniela.sforza@sbai.uniroma1.it}},\
Masahiro Yamamoto
\thanks{ 
Department of Mathematical Sciences, The University
of Tokyo, Komaba, Meguro, Tokyo 153, Japan
e-mail: {\tt myama@ms.u-tokyo.ac.jp}.
Partially supported by Grant-in-Aid for Scientific Research (S) 15H05740 of
Japan Society for the Promotion of Science}
}
\date{}
\begin{document}
\maketitle

\baselineskip 18pt

\begin{abstract}
We consider an anisotropic hyperbolic equation with memory term:
$$
\ppp_t^2 u(x,t) = \sum_{i,j=1}^n \ppp_i(a_{ij}(x)\ppp_ju)
+ \int^t_0 \sumalpha b_{\alpha}(x,t,\eta)\ppp_x^{\alpha}u(x,\eta) d\eta
+ F(x,t)
$$
for $x \in \OOO$ and $t\in (0,T)$ or $\in (-T,T)$, 
which is a model equation 
for viscoelasticity.
First we establish a Carleman estimate for this equation with 
overdetermining boundary data on a suitable lateral subboundary
$\Gamma \times (-T,T)$.
Second we apply the Carleman estimate to establish a both-sided estimate 
of $\Vert u(\cdot,0)\Vert_{H^3(\Omega)}$ by 
$\ppp_{\nu}u\vert_{\Gamma\times (0,T)}$ under the assumption that 
$\ppp_tu(\cdot,0) = 0$ and $T>0$ is sufficiently large, $\Gamma 
\subset \ppp\Omega$ satisfies some geometric condition.
Such an estimate is a kind of observability inequality and related to the 
exact controllability.  Finally we apply the Carleman estimate to an 
inverse source problem of determining a spatial varying factor 
in $F(x,t)$ and we establish a both-sided Lipschitz stability estimate.
\end{abstract}

\section{Introduction and main results}

Let $\Omega \subset \R^n$ be a bounded domain with smooth boundary 
$\ppp\Omega$.  We consider an integro-hyperbolic equation
$$
\ppp_t^2 u(x,t) = \sum_{i,j=1}^n \ppp_i(a_{ij}(x)\ppp_ju)
+ \int^t_0 \sumalpha b_{\alpha}(x,t,\eta)\ppp_x^{\alpha}u(x,\eta) d\eta
+ F(x,t), \quad x \in \OOO, 
$$
$$
\qquad \qquad \qquad \thinspace -T<t<T.       \eqno{(1.1)}
$$

We assume that 
$$\left\{\begin{array}{rl}
&a_{ij} = a_{ji} \in C^1(\ooo{\Omega}), \quad 1 \le i,j \le n, \\
&\mbox{there exists a constant $\mu_0 > 0$ such that}\\
&\sum_{i,j=1}^n a_{ij}(x)\xi_i\xi_j \ge \mu_0 \sum_{i=1}^n \xi_i^2,
\quad x \in \ooo{\Omega}, \thinspace \xi_1, ..., \xi_n \in \R.
\end{array}\right.                         \eqno{(1.2)}
$$
Here and henceforth let $\alpha = (\alpha_1, ..., \alpha_n) \in 
(\N \cup \{ 0\})^n$ be a multi-index and we set 
$\vert \alpha\vert = \alpha_1 + \cdots + \alpha_n$,
$\ppp_i = \frac{\ppp}{\ppp x_i}$, $1\le i \le n$, 
$\ppp_t = \frac{\ppp}{\ppp t}$,
$\ppp_x^{\alpha} = \ppp_1^{\alpha_1}\cdots \ppp_n^{\alpha_n}$.  

Throughout this paper, we assume
$$\left\{\begin{array}{rl}
& \nabla_{x,t}\ppp_t^kb_{\alpha} \in C(\ooo{\OOO} \times [0,T]^2),
\quad k=0,1,2 \quad \mbox{if $\vert \alpha\vert=2$}, \\
& \ppp_t^kb_{\alpha} \in L^{\infty}(\OOO\times (0,T)^2),
\quad k=0,1,2 \quad \mbox{if $\vert \alpha\vert\le 1$}, \\
\end{array}\right.            \eqno{(1.3)}
$$
if not specified.
Let $\nu = \nu(x) = (\nu_1, ..., \nu_n)$ 
be the unit outward normal vector to $\ppp\OOO$ at
$x$ and let
$$
\ppp_{\nu}u = \sum_{i,j=1}^n a_{ij}(\ppp_ju)\nu_i  \quad 
\mbox{on $\ppp\OOO$}.
$$
For concise description, we set 
$$
Au :=  \sum_{i,j=1}^n \ppp_i(a_{ij}(x)\ppp_ju), \quad
B(t,\eta)u := \sumalpha b_{\alpha}(x,t,\eta)\ppp_x^{\alpha}u(x,\eta).
$$

The equation (1.1) is a model equation for the viscoelasticity.

We point out that for some materials, the effects of memory cannot be
neglected without failing the analysis, as observed by Volterra \cite{VV}. He
embraced Boltzmann model, according to which the stress has to depend
linearly on strain history.

Our integro-differential equation (1.1) serves as a model for
describing the viscoelastic properties of  those materials whose
properties are different along several directions.
There is a huge number of papers treating viscoelastic models, as shown,
e.g., in the book Renardy, Hrusa and Nohel \cite{RHN}. 
With no pretension to be exhaustive, we cite the papers Dafermos \cite{Da} 
and Edelstein and Gurtin \cite{EG}. 
In particular, in a pioneering work \cite{Da}, Dafermos
studied an abstract formulation of our equation, giving as an application
the case of an anisotropic viscoelastic equation.

In this paper, we first establish a Carleman estimate for (1.1).
A Carleman estimate is a 
weighted $L^2$-estimate for solutions to a partial differential equation
which holds uniformly in large parameter $s>0$, and was  
derived by Carleman \cite{Ca} for proving the unique continuation property.
Second we apply the Carleman estimate to prove an estimate of initial 
value by data on suitable lateral boundary data, which is called 
an observability inequality.
Finally we discuss an inverse source problem.  More precisely, 
the external force $F$ is assumed to cause the action, but
in practice it is often that $F$ is not a priori known and so we
have to identify by available data for the sake of accurate analysis of 
the system.
We are concerned with the determination of a spatial component 
$f(x)$ of $F(x,t) := R(x,t)f(x)$ with given $R(x,t)$.
The form $R(x,t)f(x)$ is special but in applications we model the 
external force in a more special form $F(x,t) = \lambda(t)f(x)$
where $\lambda(t)$ is the time changing ratio and $f$ is the spatial
distribution of the external force. 

For the statement of the Carleman estimate, we need to introduce notations.
We set $A_0(x) = (a_{ij}(x))_{1\le i,j\le n}$ and
$$
a(x,\xi) = A_0(x)\xi\cdot \xi = \sumij a_{ij}(x)\xi_i\xi_j
                                                    \eqno{(1.4)}
$$
for all $x \in \Omega$ and $\xi = (\xi_1, ..., \xi_n) \in \R^n$.
Given functions $p(x,\xi)$ and $q(x,\xi)$, we define the Poisson 
bracket by 
$$
\{ p,q\}(x,\xi) = \sum_{j=1}^n \left(
\frac{\ppp p}{\ppp \xi_j}\frac{\ppp q}{\ppp x_j}
- \frac{\ppp p}{\ppp x_j}\frac{\ppp q}{\ppp \xi_j}\right)(x,\xi).
$$
We set 
$$
d(x) =  \vert x-x_0\vert^2, \quad x \in \R^n
$$
with fixed $x_0 \in \R^n \setminus \overline{\OOO}$.
In addition to (1.2), throughout this paper,
we assume that there exists a 
constant $\mu_1 > 0$ such that
$$
\{a, \{a, d\}\}(x,\xi) \ge \mu_1\vert A_0(x)^{-1}\xi\vert^2, \quad 
x\in\ooo{\Omega}, \thinspace \xi \in \R^n                    \eqno{(1.5)}
$$
(e.g., Bellassoued and Yamamoto \cite{BelY2}, \cite{BelY3}).
For proving a Carleman estimate, it is known that we need some condition
like (1.5), which ic called the pseudo-convexity (e.g., 
H\"ormande \cite{H}).  We refer to 
Yao \cite{Ya} which discusses anisotropic materials without intregral terms
and shows a counterexample to the observability inequality without 
such condition for the principal part. 

Next as subboundary where we take boundary data of the solution 
$u$, we define
$$
\Gamma = \{ x \in \ppp\Omega; \thinspace 
(x-x_0)\cdot \nu(x) \ge 0\}.                       \eqno{(1.6)}
$$ 
Here and henceforth $(\cdot, \cdot)$ denotes the scalar product in 
$\R^n$.

Furthermore we set 
$$
Q = \OOO \times (-T,T),
$$
$$
\psi(x,t) = \vert x-x_0\vert^2 - \beta t^2, \quad (x,t) \in \ooo{Q}
                                                    \eqno{(1.7)}
$$
and 
$$
\va(x,t) = e^{\gamma\psi(x,t)}, \quad \Phi = \Phi(\gamma) 
= \max_{(x,t)\in\ooo{Q}} \va(x,t),                   \eqno{(1.8)}
$$
where $\beta > 0$ is chosen sufficienly small for the constant $\mu_1>0$ 
in (1.5) and $\gamma>0$ is a second large paramater and chosen later.

The conditions (1.5) and (1.6) pose extra conditions for $a_{ij}$ and 
$\Gamma$ respectively and are a sufficient 
condition for the Carleman estimate below stated.

Now we introduce a cut-off function $\chi$ in $t$.
For fixed sufficiently small $\ep>0$, let $\chi \in C^{\infty}(\R)$ satisfy 
$0 \le \chi \le 1$ in $\R$ and 
$$
\chi(t) = 
\left\{\begin{array}{rl}
1, \quad & \vert t\vert \le T-2\ep, \\
0, \quad & \vert t\vert \ge T-\ep.
\end{array}\right.
                                        \eqno{(1.9)}
$$
We further set 
$$
\delta = \max_{x\in \ooo{\Omega}, 
T-2\ep \le \vert t\vert \le T-\ep} \va(x,t).         \eqno{(1.10)}
$$
Then $\va(x,t) \le \delta$ for $x\in \ooo{\Omega}$ and  
$T-2\ep \le \vert t\vert \le T-\ep$.

Now we are ready to state our first main result.
\\

{\bf Theorem 1.1.}\\
Let 
$$
F \in L^2(-T,T;H^2(\OOO)), \quad F=\ppp_{\nu}F = 0 \quad 
\mbox{on $\ppp\Omega$}.                 
$$
We set
$$
v(x,t) = \chi(t)\ppp_t^2u(x,t) - \chi(t)F(x,t), \quad x \in \Omega, 
\thinspace -T<t<T.
$$
\\
(i) There exists a constant $\gamma_0 > 0$ such that for $\gamma > \gamma_0$,
we can choose constants $s_0=s_0(\gamma)>0$ and $C>0$ such that 
$$
\int_Q ( s\gamma \va \vert\nabla_{x,t}v\vert^2 
+ s^3\gamma^3\va^3 \vert v\vert^2) \weight dxdt
$$
$$
\le C\int_Q \vert AF\vert^2 \weight dxdt 
+ C\Vert u\Vert^2_{H^1(-T,T;H^2(\OOO))} s^2\Phi^2e^{2s\delta}
+ Ce^{Cs}\Vert \ppp_{\nu}u\Vert^2_{H^2(-T,T;L^2(\Gamma))}
                                               \eqno{(1.11)}
$$
and
$$
\int_Q \left( \sumalpha s^2\va^2\vert \chi \ppp_x^{\alpha}u\vert^2
+ \vert \chi \ppp_x^{\alpha}\ppp_tu\vert^2\right) \weight dxdt
$$
$$
\le C\int_Q \vert AF\vert^2 \weight dxdt 
+ C\Vert u\Vert^2_{H^1(-T,T;H^2(\OOO))} s^2\Phi^2e^{2s\delta}
+ Ce^{Cs}\Vert \ppp_{\nu}u\Vert^2_{H^2(-T,T;L^2(\Gamma))}
                                               \eqno{(1.12)}
$$
for all $s > s_0$ and $u \in L^2(Q)$ satisfying 
$\ppp_x^{\alpha}u \in H^2(-T,T;L^2(\Omega))\cap
L^2(-T,T;H^2(\Omega))$ for all $\vert \alpha\vert \le 2$ and
$u\vert_{\ppp\OOO} = 0$.
\\
(ii) Moreover we assume
$$
F, \ppp_tF \in L^2(-T,T;H^2(\OOO)).
$$
Then 
$$
\int_Q \left( \sumalpha s^2\va^2\vert \chi \ppp_x^{\alpha}u\vert^2
+ \vert \chi \ppp_x^{\alpha}\ppp_tu\vert^2
+ \vert \chi\ppp_x^{\alpha}\ppp_t^2u\vert^2 \right) \weight dxdt
$$
$$
\le C\int_Q (\vert AF\vert^2 + \vert A\ppp_tF\vert^2)\weight dxdt 
+ C\Vert u\Vert^2_{H^2(-T,T;H^2(\OOO))} s^2\Phi^2e^{2s\delta}
+ Ce^{Cs}\Vert \ppp_{\nu}u\Vert^2_{H^3(-T,T;L^2(\Gamma))}
                                               \eqno{(1.13)}
$$
for all $s > s_0$ and $u \in H^2(Q)$ satisfying 
$\ppp_x^{\alpha}u, \ppp_tu \in H^2(Q)$ for all $\vert \alpha\vert \le 2$ and
$u\vert_{\ppp\OOO} = 0$.
\\

In (ii) of the theorem, we can rewrite (1.13) in terms of $v$, but we omit.

Inequalities (1.12) and (1.13) hold for each solution $u$ to (1.1) and both
are weighted with $e^{2s\va(x,t)}$ and uniform for sufficiently large 
$s>0$ in the sense that the constant $C>0$ is independent of all large 
$s>0$. Such an inequality is called a Carleman estimate.
The Carleman estimate is effectively applied to the unique continuation for
partial differential equations, the observability inequality and 
inverse problems.  In this paper, by Theorem 1.1 we establish the 
observability inequality (Theorem 1.2) and the Lipschitz stability in an 
inverse source problem (Theorem 1.3) for (1.1).

As for general treatments on Carleman estimates for partial differential 
equations without integral terms, we refer to H\"ormander \cite{H}, 
Isakov \cite{Is1}.  
There are many works concerning Carleman estimates for partial 
differential equations without integral terms.  
Since inverse problems are often concerned with the determination of
the principal coefficients $a_{ij}(x)$, we have to concretely 
realize the condition (1.5).  As for such concrete
Carleman estimates which give sufficient conditions for (1.5) and more 
directly applicable to inverse problems, see Amirov and Yamamoto \cite{AY}.
We refer to Baudouin, de Buhan and Ervedaza \cite{B1},
Imanuvilov \cite{Ima2}, Kha{\u\i}darov \cite{Kha}, Romanov \cite{Rom} 
which establish Carleman estimates for hyperbolic equations.
For Carleman estimates for parabolic equations, in addition to 
Isakov \cite{Is1}, \cite{Is2}, Isakov and Kim \cite{IK}, see
Fursikov and Imanuvilov \cite{FI}, Imanuvilov \cite{Ima1},
Imanuvilov, Puel and Yamamoto \cite{ImaPuY}, Yamamoto \cite{Y}.
For elliptic Carleman estimates where the right-hand side is
estimated in $H^{-1}$-space, see Imanuvilov and Puel \cite{ImaPu}.

As for isotropic hyperbolic equations 
with integral terms, Cavaterra, Lorenzi and Yamamoto \cite{Caloya}
established a Carleman estimate and applied it for proving a stability 
result for some inverse source problem.  Here the isotropic hyperbolic
equation means 
$$
a_{ij}(x) = 
\left\{ \begin{array}{rl}
p(x), \quad & i=j, \\
0,    \quad & i\ne j
\end{array}\right.
$$
in (1.1).
After \cite{Caloya}, in the case where $A$ and $B$ are isotropic Lam\'e 
operators, the 
following works discuss Carleman estimates and inverse problems:
de Buhan and Osses \cite{dBO}, Lorenzi, Messina and Romanov \cite{LoMeR}, 
Lorenzi and Romanov \cite{LoR}, Romanov and Yamamoto \cite{RY}.
However, to the best knowledge of the authors, there are no publications
on Carleman estimates for anisotropic hyperbolic equations with 
integral terms $\int^t_0 \sumalpha b_{\alpha}\ppp_x^{\alpha}u d\eta$.
For this anisotropic case, differently from the isotropic case, we need 
a lot of technicalities because $A$ and $B$ are not commutative modulo 
lower-order terms of derivatives.
\\

Now we present two applications of the Carleman estimate (Theorem 1.1).
First we derive a partial observability inequality of estimating 
one component of a pair of initial values.

Let $y \in L^2(\Omega \times (0,T))$ satisfy 
$\ppp_x^{\alpha}y \in H^2(0,T;L^2(\Omega)) \cap 
L^2(0,T;H^2(\Omega))$ with $\vert \alpha\vert \le 2$ and
$$
\left\{ \begin{array}{rl}
& \ppp_t^2y(x,t) = Ay + \int^t_0 \sumalpha b_{\alpha}(x,t,\eta)
\dalpha y(x,\eta)d\eta, \quad x \in \Omega, \thinspace 0 < t < T, \\
&y\vert_{\ppp\Omega} = 0, \quad 0 < t < T, \\
&y(x,0) = a(x), \quad \ppp_ty(x,0) = 0, \quad x\in \Omega.\\
\end{array}\right.    
                                    \eqno{(1.14)}
$$

Let $\Gamma \subset \ppp\Omega$ and $T>0$ be given.
\\
\vspace{0.2cm}

{\bf Partial observability inequality.}\\
Estimate $a(x)$ by $\ppp_{\nu}y\vert_{\Gamma\times (0,T)}$.
\\

Thanks to the integral term $\int^t_0 \sumalpha b_{\alpha}(x,t,\eta)
\dalpha y(x,\eta)d\eta$, our method requests that 
$y(\cdot,0) = 0$ or $\ppp_ty(\cdot,0) = 0$ in $\Omega$, and here we 
discuss only the case of $\ppp_ty(\cdot,0) = 0$.

The observability inequality is 
regarded as a dual problem to the exact controllability, and
for a hyperbolic type of equations without integral terms, there have been 
enormous works.  Here we refer only to Komornik \cite{Ko}, Lions \cite{L},
and Yao \cite{Ya} which discusses anisotropic hyperbolic 
equations without integral terms.  
For proving observability inequalities, the multiplier
method is commonly applied, but also a Carleman estimate is applicable for 
wider classes of partial differential equations (e.g., Kazemi and Klibanov
\cite{KK}, Klibanov and Malinsky \cite{KM}).  
As for the first application of Theorem 1.1, 
we show an observability inequality for (1.14).

{\bf Theorem 1.2.}\\
We assume that 
$\Gamma$ satisfies (1.6), and
$$
T > \frac{\max_{x\in\ooo{\Omega}} \vert x-x_0\vert}{\sqrt{\beta}}.
                                         \eqno{(1.15)}
$$
Then there exists a constant $C>0$ such that 
$$
C^{-1}\Vert \ppp_{\nu}y\Vert_{H^2(0,T;L^2(\ppp\OOO))}
\le \Vert y(\cdot,0)\Vert_{H^3(\Omega)} 
\le C\Vert \ppp_{\nu}y\Vert_{H^2(0,T;L^2(\Gamma))}.
$$
for each solution $y$ to (1.14) with $y(\cdot,0) \in H^3(\OOO) 
\cap H^1_0(\OOO)$.
\\
\vspace{0.3cm}

In this theorem, we can replace (1.3) by weaker condition
$$\left\{\begin{array}{rl}
&\ppp_t^k\nabla b_{\alpha} \in C(\ooo{\Omega}\times [0,T]^2), \quad

k=0,1,2, \thinspace \vert \alpha\vert = 2, \\
&\ppp_t^kb_{\alpha} \in C(\ooo{\Omega}\times [0,T]^2), \quad
k=0,1,2, \thinspace \vert \alpha\vert \le 1. \\
\end{array}\right.
$$
 
By the finiteness of the propagation speed, we need to assume (1.15), and
also a geometric condition (1.6) on the observation subboundary $\Gamma$ is 
assumed.
This is the same for the inverse source problem stated below.

Finally we discuss an inverse source problem.  That is, we consider
$$
\left\{ \begin{array}{rl}
& \ppp_t^2u = Au + \int^t_0 \sumalpha b_{\alpha}(x,t,\eta)\dalpha u(x,\eta)
d\eta + R(x,t)f(x), \quad x \in \Omega, \thinspace 0 < t < T, \\
& u\vert_{\ppp\OOO} = 0, \quad 0 < t <  T,\\
& u(x,t) = \ppp_tu(x,0) = 0, \quad x \in \OOO.
\end{array}\right.                            \eqno{(1.16)}
$$
Here we assume
$$\left\{\begin{array}{rl}
& R \in H^3(0,T;W^{2,\infty}(\OOO)), \\
& \vert R(x,0) \vert \ne 0, \quad x \in \ooo{\OOO}.
\end{array}\right.                            \eqno{(1.17)}
$$
Let $\Gamma \subset \ppp\OOO$ be given and $T>0$ be fixed.
Then we discuss\\

{\bf Inverse source problem.}\\
Determine $f(x)$, $x \in\OOO$ from $\ppp_{\nu}u\vert_{\Gamma\times (0,T)}$.
\\

As the stability for the inverse problem, we prove

{\bf Theorem 1.3.}\\
We assume (1.6) and (1.15).
Then there exists a constant $C>0$ such that 
$$
C^{-1}\Vert \ppp_{\nu}u\Vert_{H^3(0,T;L^2(\ppp\OOO))}
\le \Vert f\Vert_{H^2(\OOO)} \le C\Vert \ppp_{\nu}u\Vert
_{H^3(0,T;L^2(\Gamma))}
                                      \eqno{(1.18)}
$$
for each $f \in H^2_0(\OOO)$.
\\

The second inequality in (1.18) asserts the Lipschitz stability 
for our inverse problem.   The first inequality means that
our estimate is the best possible estimate for the inverse source 
problem.

Our argument for the inverse problem is 
based on Bukhgeim and Klibanov \cite{BK}, which relies on a Carleman estimate.
Klibanov \cite{Kli2} corresponds to the full version of \cite{BK}.
Since Bukhgeim and Klibanov \cite{BK}, their methodology has been 
developed for various equations and we can refer to many papers on 
inverse problems of determining spatially varying coefficients and 
components of source terms.  As a partial list of references on inverse 
problems for hyperbolic and parabolic equations by Carleman estimates, 
we refer to Baudouin and Yamamoto \cite{BY}, Bellassoued \cite{Bel1},
\cite{Bel2}, Bellassoued and Yamamoto \cite{BelY1},
Benabdallah, Cristofol, Gaitan and Yamamoto \cite{BCGY},
Cristofol, Gaitan and Ramoul \cite{CGR}, Imanuvilov and Yamamoto
\cite{IY1} - \cite{IY4}, Klibanov \cite{Kli1}, \cite{Kli2},
Klibanov and Yamamoto \cite{KY}, Yamamoto \cite{Y},
Yuan and Yamamoto \cite{YY2}.

As for similar inverse problems for the Navier-Stokes equations, 
see Bellassoued, Imanuvilov
and Yamamoto \cite{BelIY}, Choulli, Imanuvilov, Puel and Yamamoto
\cite{CIPY}, Fan, Di Cristo, Jiang and Nakamura \cite{FDJN}, 
Fan, Jiang and Nakamura \cite{FJN}.
Gaitan and Ouzzane \cite{GO}, and G\"olgeleyen and Yamamoto 
\cite{GY} discuss inverse problems for transport equations by 
Carleman estimate, and 
Imanuvilov, Isakov and Yamamoto \cite{ImIsY}, Imanuvilov and Yamamoto
\cite{IY5} - \cite{IY7} discuss Carleman estimates and inverse problems 
for non-stationary isotropic Lam\'e systems, which are related to our 
equation for the viscoelasticity.
See Yuan and Yamamoto \cite{YY1} about a Carleman estimate and inverse
problems for a plate equation.
As related books on Carleman estimates and inverse problems,
see Beilina and Klibanov \cite{Bekli}, Klibanov and Timonov \cite{KT},
Lavrent'ev, Romanov and Shishat$\cdot$ski\u\i \cite{LRS}.
 
This paper is composed of five sections.  In section 2, we prove 
Theorem 1.1 and section 3 is devoted to providing fundamental 
energy estimates.  In sectios 4 and 5, we prove Theorems 1.2 and 1.3
respectively.

\section{Proof of Theorem 1.1}

The proof is combination of hyperbolic and elliptic Carleman 
estimates (Lemmata 2.1 and 2.2) with another key lemma 
(Lemma 2.3) which can incorporate the integral term in (1.1).
We divide the proof into five steps.

{\bf First Step.}\\
Henceforth we set 
$$
\Sigma = \Gamma\times (-T,T).
$$
Under the assumption (1.7), a Carleman estimate for hyperbolic equations
is known (e.g., Bellassoued and Yamamoto \cite{BelY2}, \cite{BelY3}).
\\
\vspace{0.2cm}
{\bf Lemma 2.1.}\\
There exists a constant $\gamma_0 > 0$ such that for $\gamma > \gamma_0$,
we can choose constants $s_0=s_0(\gamma) > 0$ and $C = C(\gamma) >0$ 
such that 
$$
\int_Q ( s\va\gamma\vert \nabla_{x,t}u\vert^2 
+ s^3\va^3\gamma^3\vert u\vert^2) \weight dxdt
$$
$$
\le C\int_Q \vert (\ppp_t^2 - A)u\vert^2 \weight dxdt 
+ Ce^{Cs}\Vert \ppp_{\nu}u\Vert^2_{L^2(\Sigma)}
                                               \eqno{(2.1)}
$$
for all $s > s_0$ and $u \in H^2(Q)$ satisfying 
$u\vert_{\ppp\OOO} = 0$.
\\

Moreover we have a Carleman estimate for the elliptic operator
$A$ without the extra conditions on $a_{ij}$.\\
\vspace{0.2cm}
{\bf Lemma 2.2.}\\
Let $p \in \R$ be given.  There exists a constant $\gamma_0 > 0$ such that for 
$\gamma > \gamma_0$,
we can choose constants $s_0>0$ and $C>0$ such that 
$$
\int_Q \left( s^p\va^p\sum_{\vert\alpha\vert=2}
\vert \ppp_x^{\alpha}y\vert^2 + s^{p+2}\gamma^2\va^{p+2}\vert \nabla y\vert^2
+ s^{p+4}\gamma^4\va^{p+4}\vert y\vert^2\right) \weight dxdt
$$
$$
\le C\int_Q s^{p+1}\va^{p+1}\vert Ay\vert^2 \weight dxdt 
+ Ce^{Cs}\Vert \ppp_{\nu}y\Vert^2_{L^2(\Sigma)}
                                               \eqno{(2.2)}
$$
for all $s > s_0$ and $y \in L^2(-T,T;H^2(\OOO) \cap H^1_0(\OOO))$.
\\

In the case of $p=-1$, Lemma 2.2 is classical and we refer to
Lemma 7.1 in Bellassoued and Yamamoto \cite{BelY3} for example.
For completeness, we give the proof of Lemma 2.2 for arbitrary $p \in \R$ on 
the basis of the case of $p=-1$ in Appendix.
\\
\vspace{0.2cm}
{\bf Second Step.}

For gaining compact supports in time for functions under 
consideration, we use the cut-off function.  That is, we 
recall that we choose a cut-off function $\chi\in 
C^{\infty}(\R)$ such that $0 \le \chi\le 1$ and 
$$
\chi(t) =
\left\{\begin{array}{rl}
1, \quad &\vert t\vert \le T - 2\ep, \\
0, \quad &\vert t\vert \ge T-\ep.
\end{array}\right.                             \eqno{(2.3)}
$$
In the succeeding arguments, we notice that all the terms with derivatives
of $\chi$ can be regarded as of minor orders with respect to 
the large parameter $s$.

For treating an integral, it is essential to introduce a new function
$$
v(x,t) := \chi(t)Au(x,t) + \chi(t)\int^t_0 B(t,\eta)u(x,\eta) d\eta.
                                               \eqno{(2.4)}
$$
Then
$$
v(x,t) = \chi(t)\ppp_t^2u(x,t) - \chi(t)F(x,t), \quad (x,t) \in Q.
                                                               \eqno{(2.5)}
$$
We set 
$$
\www{B}(t) = B(t,t) \quad \mbox{in $Q$}.   
$$
Then we have
\begin{align*}
& \ppp_tv(x,t) = \chi A\ppp_tu + \chi \www{B}(t)u\\
+& \chi\int^t_0 \ppp_tB(t,\eta)u d\eta
+ \chi'Au + \chi'\int^t_0 B(t,\eta)u d\eta,
\end{align*}
and so
$$
\ppp_t^2 v(x,t) = \chi\left(A\ppp_t^2u + \www{B}\ppp_tu
+ (\ppp_t\www{B})u + (\ppp_tB)(t,t)u
+ \int^t_0 \ppp_t^2B(t,\eta)u(x,\eta)d\eta\right) \eqno{(2.6)}
$$
\begin{align*}
+ & 2\chi'(t)\left( A\ppp_tu + \www{B}u 
+ \int^t_0 (\ppp_tB(t,\eta))u d\eta \right)
+ \chi''(t)\left( Au + \int^t_0 B(t,\eta)u d\eta\right)\\
=: & \chi\left(A\ppp_t^2u + \www{B}\ppp_tu
+ (\ppp_t\www{B})u + (\ppp_tB)(t,t)u
+ \int^t_0 \ppp_t^2B(t,\eta)u(x,\eta)d\eta\right) \\
+& S(x,t),
\end{align*}
Then $S(x,t)$ satisfies  
$$
\vert S(x,t)\vert \le C(\vert \chi'(t)\vert + \vert \chi''(t)\vert)
\sumalpha \sum_{k=0}^1 \left(
\vert \ppp_x^{\alpha}\ppp_t^k u(x,t)\vert
+ \left\vert\int^t_0 \vert \ppp_x^{\alpha}u(x,\eta)\vert d\eta\right\vert
\right)\quad \mbox{in $Q$}.             \eqno{(2.7)}
$$
By (2.3), we note that $\vert S(x,t)\vert \ne 0$ only if
$T - 2\ep \le \vert t\vert \le T-\ep$.
Moreover, since
$$
Av(x,t) = \chi(t)A\left( Au(x,t) + \int^t_0 B(x,t,\eta)u(x,\eta) d\eta\right),
$$
we obtain
\begin{align*}
& \ppp_t^2v - Av\\
=& \chi A\left( \ppp_t^2u - Au(x,t) - \int^t_0 B(x,t,\eta)u(x,\eta) d\eta
\right)\\
+ & \chi\left(\www{B}\ppp_tu + (\ppp_t\www{B})u + (\ppp_tB)(t,t)u
+ \int^t_0 (\ppp_t^2B)(t,\eta)u d\eta\right)\\
+& S(x,t).
\end{align*}
By $F\vert_{\ppp\OOO} = 0$ and (2.5), we see that $v\vert_{\ppp\OOO} 
= 0$.  Thus
$$
\left\{ \begin{array}{rl}
& \ppp_t^2v - Av = \chi AF + \chi J + S \quad\mbox{in $Q$},\\
& \ppp_t^jv(\cdot, \pm T) = 0 \quad \mbox{in $\OOO$, $j=0,1$},\\
& v\vert_{\ppp\OOO} = 0,
\end{array}\right.
                                           \eqno{(2.8)}
$$
where
$$
J(x,t) = \sumalpha \sum_{k=0}^1 c_{\alpha,k}\ppp_x^{\alpha}
\ppp_t^ku(x,t) + \int^t_0 \sumalpha
\www{c_{\alpha}}(x,t,\eta)\ppp_x^{\alpha}u(x,\eta) d\eta,
\quad (x,t) \in Q,                      \eqno{(2.9)}
$$
and
$c_{\alpha,k}\in W^{1,\infty}(-T,T;L^2(\OOO))$,
$\www{c}_{\alpha} \in W^{1,\infty}(\OOO \times (-T,T)^2)$.

Applying Lemma 2.1 to (2.8), we have
$$
\int_Q (s\gamma\va \vert \nabla_{x,t}v\vert^2 + s^3\va^3\gamma^3
\vert v\vert^2) \weight dxdt 
\le C\int_Q \chi^2 \vert AF\vert^2 \weight dxdt  \eqno{(2.10)}
$$
\begin{align*}
+ & C\int_Q \chi^2 \vert J\vert^2 \weight dxdt
+ C\int_Q \vert S\vert^2 \weight dxdt 
+ Ce^{Cs}\Vert \ppp_{\nu}v\Vert^2_{L^2(\Sigma)}\\
\le& C\int_Q \chi^2\vert AF\vert^2 \weight dxdt\\
+ & C\int_Q \chi^2 \sumalpha \sum_{k=0}^1 \vert \ppp_x^{\alpha}
\ppp_t^ku\vert^2 \weight dxdt
+ C\int_Q \chi^2 \sumalpha \left\vert
\int^t_0 \vert \ppp_x^{\alpha}u(x,\eta)\vert^2 d\eta\right\vert
\weight dxdt\\
+ & C\int_Q \vert S\vert^2 \weight dxdt 
+ Ce^{Cs}\Vert \ppp_{\nu}v\Vert^2_{L^2(\Sigma)}
\end{align*}
for $s > s_0$.
\\
\vspace{0.2cm}
{\bf Third Step.}\\
For estimating the integral term in (1.1) with the weight 
$e^{2s\va}$, we need to prove
\\
{\bf Lemma 2.3.}\\
Let $q\ge 0$.  Then
\begin{align*}
& \int_Q \chi^2(t)(s\va)^q \left\vert \int^t_0 
\vert w(x,\eta)\vert d\eta \right\vert^2 \weight dxdt\\
\le & C\int_Q s^{q-1}\gamma^{-1}\va^{q-1}\chi^2 \vert w(x,t)\vert^2
\weight dxdt\\
+ & \int_Q s^{q-1}\gamma^{-1}\va^{q-1}\vert \ppp_t(\chi^2)\vert 
\weight \left\vert \int^t_0 \vert w(x,\eta)\vert^2 d\eta \right\vert
dxdt.
\end{align*}
\\
This type of inequality is essential for applications of  
Carleman estimates to inverse problems (Bukhgeim and Klibanov 
\cite{BK}, Klibanov \cite{Kli2}) and the inequality not 
involving the cut-off function $\chi$, is proved in 
\cite{Kli2}, \cite{KT}.
\\
{\bf Proof.}\\
It suffices to prove for $t\ge 0$, because the proof for $t \le 0$ is 
similar.
By the Cauchy-Schwarz inequality, we have
$$
\int_Q \chi^2(t)(s\va)^q \left\vert \int^t_0 
\vert w(x,\eta)\vert d\eta \right\vert^2 \weight dtdx
                                           \eqno{(2.11)}
$$
$$
\le \int_{\OOO}\int^T_0 \chi^2(t) \left(
\int^t_0 \vert w(x,\eta)\vert^2 d\eta \right) s^q\va^qte^{2s\va}dtdx.
$$
Noting that
$$
s^q\va^qte^{2s\va(x,t)}
= - \frac{\ppp_t(e^{2s\va(x,t)})}{4\gamma\beta}(s\va)^{q-1},
$$
by integration by parts, we obtain
$$
\int_{\OOO}\int^T_0 \chi^2(t) \left(
\int^t_0 \vert w(x,\eta)\vert^2 d\eta \right) 
s^q\va^qt e^{2s\va}dxdt 
                                     \eqno{(2.12)}
$$
\begin{align*}
= & \int_{\OOO}\int^T_0  -\frac{\ppp_t(e^{2s\va(x,t)})}{4\gamma\beta}
(s\va)^{q-1}\chi^2(t)\left(\int^t_0 w^2 d\eta\right) dtdx\\
=& \int_{\OOO} \left[ \frac{e^{2s\va(x,t)}}{4\gamma\beta}(s\va)^{q-1}
\chi^2(t)\int^t_0 w^2 d\eta \right]^{t=0}_{t=T} dx\\
+& \int_{\OOO}\int^T_0 \frac{e^{2s\va(x,t)}}{4\gamma\beta}
\ppp_t((s\va)^{q-1}\chi^2(t))\left( \int^t_0 w^2d\eta\right)dtdx\\
+& \int_{\OOO}\int^T_0 \frac{e^{2s\va(x,t)}}{4\gamma\beta}(s\va)^{q-1}
\chi^2(t)w^2 dtdx \\
=& \int_{\OOO} -\frac{q-1}{2}(s\va)^{q-1}t\chi^2(t) \left(
\int^t_0 w^2 d\eta\right) \weight dtdx
+ \int_{\OOO}\int^T_0 \frac{(s\va)^{q-1}\ppp_t(\chi^2)}
{4\gamma\beta}
\weight \left( \int^t_0 w^2 d\eta\right) dtdx\\
+& \int_{\OOO}\int^T_0 \frac{e^{2s\va}}{4\gamma\beta}(s\va)^{q-1}
\chi^2(t)w^2 dtdx.
\end{align*}
Here we used 
$$
\int_{\OOO} \left[ \frac{e^{2s\va(x,t)}}{4\gamma\beta}(s\va)^{q-1}
\chi^2(t)\int^t_0 w^2 d\eta \right]^{t=0}_{t=T} dx
= 0
$$
by (2.3).

Therefore we can shift the first term on the right-hand side into the 
left-hand side, we have
\begin{align*}
&\int_{\OOO}\int^T_0 t s^q\va^q \left(1 - \frac{\vert q-1\vert}{2}
\frac{1}{s\va} \right) \chi^2(t) \left(
\int^t_0 \vert w(x,\eta)\vert^2 d\eta \right) e^{2s\va}dtdx\\
\le &\int_{\OOO}\int^T_0 t s^q\va^q \left(1 + \frac{q-1}{2}
\frac{1}{s\va} \right) \chi^2(t) \left(
\int^t_0 \vert w(x,\eta)\vert^2 d\eta \right) e^{2s\va}dtdx\\
= & \int_{\OOO}\int^T_0 \frac{(s\va)^{q-1}\ppp_t(\chi^2)}
{4\gamma\beta}
\weight \left( \int^t_0 w^2 d\eta\right) dtdx\\
+& \int_{\OOO}\int^T_0 \frac{e^{2s\va}}{4\gamma\beta}(s\va)^{q-1}
\chi^2(t)w^2 dtdx.
\end{align*}

Choosing $\gamma > 0$ and $s>0$ sufficiently large and noting 
that $\va = e^{\gamma\psi}$ and $\psi \ge 0$ in $Q$, we can 
obtain $1 - \frac{\vert q-1\vert}{2}\frac{1}{s\va} \ge \frac{1}{2}$.
Therefore
\begin{align*}
&\int_{\OOO}\int^T_0 t s^q\va^q \chi^2(t) \left(
\int^t_0 \vert w(x,\eta)\vert^2 d\eta \right) e^{2s\va}dtdx\\
\le& C\int_{\OOO}\int^T_0 \frac{(s\va)^{q-1}}{\gamma}\vert\ppp_t(\chi^2)\vert
\weight \left( \int^t_0 w^2 d\eta\right) dtdx\\
+& C\int_{\OOO}\int^T_0 \frac{e^{2s\va}}{\gamma}(s\va)^{q-1}
\chi^2(t)w^2 dtdx.
\end{align*}
Substituting this into (2.12), by (2.11) we can complete the proof of 
Lemma 2.3.
\\
\vspace{0.2cm}
{\bf Fourth Step.}\\
Henceforth $\mu(t)$ generically denotes functions in $L^{\infty}(\R)$ such 
that 
$$
\mu(t)\ge 0, \quad 
\mu(t) \ne 0 \quad \mbox{only if $T - 2\ep \le \vert t\vert \le T-\ep$}.
$$
Henceforth we denote $\chi'(t) = \frac{d\chi}{dt}(t)$,
$\chi''(t) = \frac{d^2\chi}{dt^2}(t)$.
We note that $\vert \ppp_t(\chi^2)\vert$, $\vert \chi'(t)\vert^2$, 
$\vert \chi''(t)\vert^2$ can be
replaced by $\mu(t)$ in the following estimation.

Applying Lemma 2.3 to the third term on the right-hand side of (2.10) and
noting that $\chi\ppp_x^{\alpha}\ppp_tu 
= \ppp_x^{\alpha}\ppp_t(\chi u) - \chi'\ppp_x^{\alpha}u$,
we have
\begin{align*}
& \int_Q (s\gamma\va\vert \nabla_{x,t}v\vert^2 + s^3\gamma^3\va^3\vert v
\vert^2) \weight dxdt
\le C\int_Q \chi^2\vert AF\vert^2 \weight dxdt\\
+& C\int_Q \sumalpha \vert \ppp_x^{\alpha}\ppp_t(\chi u)
- \chi'\ppp_x^{\alpha}u\vert^2 \weight dxdt
+ C\int_Q \sumalpha \vert \ppp_x^{\alpha}(\chi u)\vert^2 \weight dxdt\\
+ & C\sumalpha \left( \int_Q s^{-1}\gamma^{-1}\va^{-1}
\vert \dalpha (\chi u)\vert^2 \weight dxdt
+ \int_Q s^{-1}\gamma^{-1}\va^{-1}\vert \ppp_t(\chi^2)\vert \weight
\left\vert \int^t_0 \vert \dalpha u\vert^2 d\eta \right\vert dxdt
\right)\\
+& C\int_Q (\vert \chi'\vert^2+\vert \chi''\vert^2)\weight
\sumalpha \sum_{k=0}^1 
\left( \vert \dalpha\ppp_t^ku\vert^2 
+ \left\vert \int^t_0 \vert \dalpha\ppp_t^ku \vert^2d\eta\right\vert
\right) 
+ Ce^{Cs}\Vert \ppp_{\nu}v\Vert^2_{L^2(\Sigma)}.
\end{align*}
We set
$$
U_1(x,t) = \sumalpha \sum_{k=0}^1 
\left( \vert \dalpha\ppp_t^ku(x,t)\vert^2 
+ \left\vert \int^t_0 \vert \dalpha\ppp_t^ku (x,\eta)\vert^2 
d\eta\right\vert \right).
$$
Then 
\begin{align*}
& \int_Q (s\gamma\va\vert \nabla_{x,t}v\vert^2 + s^3\gamma^3\va^3\vert v
\vert^2) \weight dxdt
\le C\int_Q \chi^2\vert AF\vert^2 \weight dxdt\\
+& C\int_Q \sumalpha \vert \ppp_x^{\alpha}\ppp_t(\chi u)\vert^2
\weight dxdt\\
+ &C\int_Q \sumalpha \vert \ppp_x^{\alpha}(\chi u)\vert^2 \weight dxdt\\
+ & C\int_Q (\vert \chi'\vert^2+\vert \chi''\vert^2 
+ \vert \ppp_t(\chi^2)\vert)U_1(x,t)\weight dxdt
+ Ce^{Cs}\Vert \ppp_{\nu}v\Vert^2_{L^2(\Sigma)}\\
+& C\sumalpha \int_Q s^{-1}\gamma^{-1}\va^{-1}
\vert \ppp_t(\chi^2)\vert \weight 
\left\vert \int^t_0 \vert \dalpha u\vert^2 d\eta\right\vert dxdt\\
+& C\int_Q (\vert \chi'\vert^2 + \vert \chi''\vert^2)U_1\weight dxdt
+ Ce^{Cs}\Vert \ppp_{\nu}v\Vert^2_{L^2(\Sigma)}.
\end{align*}
Therefore
$$
\int_Q (s\gamma\va\vert \nabla_{x,t}v\vert^2 + s^3\gamma^3\va^3\vert v
\vert^2) \weight dxdt
\le C\int_Q \chi^2\vert AF\vert^2 \weight dxdt   \eqno{(2.13)}
$$
\begin{align*}
+ & C\int_Q \sumalpha \vert \ppp_x^{\alpha}(\chi \ppp_tu)\vert^2
\weight dxdt
+ C\int_Q \sumalpha \vert \ppp_x^{\alpha}(\chi u)\vert^2 \weight dxdt\\
+ & C\int_Q \mu(t)U_1\weight dxdt 
+ Ce^{Cs}\Vert \ppp_{\nu}v\Vert^2_{L^2(\Sigma)}
\end{align*}
for $s > s_0$.
\\

We will estimate the second and the third terms on the right-hand side
of (2.13).  By (2.4), we have
$$
A(\chi(t)u(x,t)) = v(x,t) - \chi(t)\int^t_0 B(t,\eta)u(x,\eta) d\eta
$$
in $Q$, and so we apply (2.2) with $p=2$, we obtain
\begin{align*}
& \int_Q s^2\va^2 \sumalpha \vert \dalpha (\chi u)\vert^2 \weight dxdt\\
\le& C\int_Q s^3\va^3\vert v\vert^2 \weight dxdt
+ C\int_Q s^3\va^3 \left\vert \chi(t)\int^t_0 B(t,\eta)u(x,\eta) d\eta
\right\vert^2 \weight dxdt\\
+ & Ce^{Cs}\Vert \ppp_{\nu}u\Vert^2_{L^2(\Sigma)}.
\end{align*}

Next we apply Lemma 2.3 with $q =2$ to the second term on the 
right-hand side and, similarly to (2.13), we obtain
$$
\int_Q s^3\va^3\chi^2(t) \left\vert \int^t_0 \vert B(t,\eta)u\vert d\eta
\right\vert^2 \weight dxdt              \eqno{(2.14)}
$$
\begin{align*}
\le& C\int_Q s^3\va^3\chi^2(t)\left\vert \int^t_0 \sumalpha
\vert \dalpha u\vert^2 d\eta \right\vert \weight dxdt
\le C\int_Q s^3\va^3\chi^2(t)\sumalpha \left\vert \int^t_0
\vert \dalpha u\vert^2 d\eta \right\vert \weight dxdt\\
\le & C\sumalpha \int_Q s^2\gamma^{-1}\va^2 \chi^2(t)
\vert \dalpha u\vert^2 \weight dxdt
+ C\sumalpha \int_Q s^2\gamma^{-1}\va^2 \vert \ppp_t(\chi^2)\vert
\weight \left\vert \int^t_0 \vert \dalpha u\vert^2 d\eta \right\vert dxdt\\
\le & C\sumalpha \int_Q s^2\gamma^{-1}\va^2 \chi^2(t)
\vert \dalpha u\vert^2 \weight dxdt
+ C\int_Q s^2\gamma^{-1}\va^2 \vert \ppp_t(\chi^2)\vert
\weight U_1(x,t) dxdt.
\end{align*}
Therefore
\begin{align*}\
& \int_Q s^2\va^2 \sumalpha \vert \dalpha (\chi u)\vert^2
\weight dxdt\\
\le& C\int_Q s^3\va^3 \vert v\vert^2 \weight dxdt
+ C\sumalpha \int_Q s^2\gamma^{-1}\va^2\vert \dalpha (\chi u)\vert^2
\weight dxdt\\
+& C\int_Q s^2\gamma^{-1}\va^2 \mu(t)U_1 \weight dxdt
+ Ce^{Cs}\Vert \ppp_{\nu}u\Vert^2_{L^2(\Sigma)}.
\end{align*}

Choosing $\gamma > 0$ sufficiently large, we can absorb the second term 
on the right-hand side into the left-hand side.
Moreover we choose $s_0(\gamma) > 0$ sufficiently larger such that
$s_0(\gamma)\va(x,t) \ge s_0(\gamma)e^{-\gamma \max_{\overline{Q}}
\vert \psi(x,t)\vert} \ge 1$, and by $s\va \ge 1$ we obtain
$$
\int_Q \sumalpha \vert \dalpha (\chi u)\vert^2 \weight dxdt
\le \int_Q s^2\va^2 \sumalpha \vert \dalpha (\chi u)\vert^2
\weight dxdt                       \eqno{(2.15)}
$$
\begin{align*}
\le& C\int_Q s^3\va^3 \vert v\vert^2 \weight dxdt
+ C\int_Q s^2\gamma^{-1}\va^2 \mu(t) U_1\weight dxdt
+ Ce^{Cs}\Vert \ppp_{\nu}u\Vert^2_{L^2(\Sigma)}
\end{align*}
for $s > s_0$.

Now we estimate the third term on the right-hand side of (2.13).  
By (2.4) we have
\begin{align*}
& A(\chi \ppp_tu) = \ppp_tv - \chi\www{B}u 
- \chi\int^t_0 (\ppp_tB)(t,\eta)u d\eta\\
-& \chi'Au - \chi'\int^t_0 B(t,\eta)u d\eta.
\end{align*}
Apply (2.2) with $p=0$, and we obtain
\begin{align*}
& \int_Q \sumalpha \vert \dalpha (\chi \ppp_tu)\vert^2 \weight dxdt\\
\le & C\int_Q s\va \vert \ppp_tv\vert^2 \weight dxdt
+ C\int_Q s\va \chi^2 \vert \www{B}u\vert^2 \weight dxdt\\
+& C\int_Q s\va\chi^2 \left\vert \int^t_0
(\ppp_tB)(t,\eta)u d\eta \right\vert^2\weight dxdt
+ C\int_Q \vert \chi'\vert^2 s\va \vert Au\vert^2 \weight dxdt\\
+& C\int_Q \vert \chi'\vert^2s\va
\left\vert \int^t_0 Bu d\eta\right\vert^2 \weight dxdt
+ Ce^{Cs}\Vert \ppp_{\nu}(\chi\ppp_tu)\Vert^2_{L^2(\Sigma)}.
\end{align*}
By an argument similar to (2.14) for the third term on the 
right-hand side, we have
\begin{align*}\
& \int_Q \sumalpha \vert \dalpha (\chi \ppp_tu)\vert^2
\weight dxdt\\
\le& C\int_Q s\va \vert \ppp_tv\vert^2 \weight dxdt
+ C\int_Q s\va\chi^2 \vert \www{B}u\vert^2 \weight dxdt\\
+& C\left(\sumalpha \int_Q \gamma^{-1}\chi^2(t)\vert \dalpha u\vert^2
\weight dxdt
+ \sumalpha \int_Q \gamma^{-1}\vert \ppp_t(\chi^2)\vert \
\left\vert \int^t_0 \vert \dalpha u\vert^2 d\eta \right\vert
\weight dxdt\right) \\
+& C\int_Q \vert \chi'\vert^2 s\va \vert Au\vert^2 \weight dxdt
+ C\int_Q \vert \chi'\vert^2 s\va \left\vert 
\int^t_0 \vert B(t,\eta)u \vert d\eta \right\vert^2 \weight dxdt
+ Ce^{Cs}\Vert \ppp_{\nu}(\chi\ppp_tu)\Vert^2_{L^2(\Sigma)}\\
\le& C\int_Q s\va \vert \ppp_tv\vert^2 \weight dxdt
+ C\sumalpha \int_Q s\va \vert \dalpha (\chi u)\vert^2 \weight dxdt
+ C\sumalpha \int_Q \gamma^{-1} \vert \dalpha (\chi u)\vert^2 \weight dxdt\\
+ & C\int_Q s\va \mu(t) U_1 \weight dxdt
+ Ce^{Cs}\Vert \ppp_{\nu}\ppp_tu\Vert^2_{L^2(\Sigma)}.
\end{align*}
Applying (2.15), we obtain
$$
\int_Q \sumalpha \vert \dalpha (\chi \ppp_tu)\vert^2
\weight dxdt \le C\int_Q (s\va \vert \ppp_tv\vert^2 + s^2\va^2\vert v\vert^2)
\weight dxdt
                                 \eqno{(2.16)}
$$
\begin{align*}
+ &C\int_Q U_1(x,t)s^2\gamma^{-1}\va^2 \mu(t) \weight dxdt
+ Ce^{Cs}(\Vert \ppp_{\nu}u\Vert^2_{L^2(\Sigma)}
+ \Vert \ppp_{\nu}\ppp_tu\Vert^2_{L^2(\Sigma)})
\end{align*}
for $s > s_0$.

Applying (2.15) and (2.16) to (2.13), we obtain
\begin{align*}
& \int_Q (s\gamma\va \vert \nabla_{x,t}v\vert^2 
+ s^3\gamma^3\va^3 \vert v\vert^2) \weight dxdt\\
\le &C\int_Q \vert AF\vert^2 \weight dxdt
+ C\int_Q (s\va \vert \nabla_{x,t}v\vert^2 
+ s^3\va^3 \vert v\vert^2) \weight dxdt\\
+ & C\int_Q s^2\gamma^{-1}\va^2 \mu(t) U_1 \weight dxdt
+ Ce^{Cs}(\Vert \ppp_{\nu}u\Vert^2_{L^2(\Sigma)}
+ \Vert \ppp_{\nu}\ppp_tu\Vert^2_{L^2(\Sigma)}
+ \Vert \ppp_{\nu}v\Vert^2_{L^2(\Sigma)}).
\end{align*}
Choosing $\gamma > 0$ sufficiently large, we can absorb the second
term on the right-hand side into the left-hand side.
Therefore
$$
\int_Q (s\gamma\va \vert \nabla_{x,t}v\vert^2 
+ s^3\gamma^3\va^3 \vert v\vert^2) \weight dxdt
                                         \eqno{(2.17)}
$$
\begin{align*}
\le &C\int_Q \vert AF\vert^2 \weight dxdt
+ C\int_Q s^2\va^2 \mu(t) U_1 \weight dxdt\\
+ &Ce^{Cs}(\Vert \ppp_{\nu}u\Vert^2_{L^2(\Sigma)}
+ \Vert \ppp_{\nu}\ppp_tu\Vert^2_{L^2(\Sigma)}
+ \Vert \ppp_{\nu}v\Vert^2_{L^2(\Sigma)})
\end{align*}
for $s > s_0$.
By the definition of $U_1$, the estimate (2.17) proves (1.11).

Next we will prove (1.12). 
The addition of (2.15) and (2.16) yields
$$
\int_Q (s^2\va^2 \vert \dalpha(\chi u)\vert^2
+ \vert \dalpha(\chi \ppp_tu)\vert^2) \weight dxdt
                                      \eqno{(2.18)}
$$
\begin{align*}
\le &C\int_Q (s\va\vert \ppp_tv\vert^2 + s^3\va^3\vert v\vert^2)\weight
dxdt
+ C\int_Q s^2\va^2 \mu(t) U_1 \weight dxdt\\
+ & Ce^{Cs}(\Vert \ppp_{\nu}u\Vert^2_{L^2(\Sigma)}
+ \Vert \ppp_{\nu}\ppp_tu\Vert^2_{L^2(\Sigma)}).
\end{align*}
Applying (2.17) in (2.18), we obtain
\begin{align*}
& \int_Q \sumalpha (
s^2\va^2 \vert \dalpha (\chi u)\vert^2 
+ \vert \dalpha(\chi \ppp_tu)\vert^2)\weight dxdt\\
\le & C\int_Q \vert AF\vert^2 \weight dxdt\\
+& C\int_Q s^2\va^2 \mu(t) U_1 \weight dxdt
+ Ce^{Cs}(\Vert \ppp_{\nu}u\Vert^2_{L^2(\Sigma)}
+ \Vert \ppp_{\nu}\ppp_tu\Vert^2_{L^2(\Sigma)}
+ \Vert \ppp_{\nu}v\Vert^2_{L^2(\Sigma)}).
\end{align*}
By (2.4) and $\ppp_{\nu}F\vert_{\ppp\Omega}=0$, we have
$\ppp_{\nu}v = \chi\ppp_{\nu}\ppp_t^2u - \chi\ppp_{\nu}F
= \chi\ppp_{\nu}\ppp_t^2u$ on $\ppp\Omega$.  
Consequently the estimate (1.12), and the proof of Theorem 1.1 (i)
is completed.
\\

{\bf Fifth Step.}\\
We will prove Theorem 1.1 (ii).
By (2.4), we have
$$\left\{ \begin{array}{rl}
& A(\chi \ppp_t^2u) = \ppp_t^2v \\
-& \chi\left( \www{B}\ppp_tu + (\ppp_t\www{B})(t)u
+ (\ppp_tB)(t,t)u + \int^t_0 (\ppp_t^2B)u d\eta\right)\\
-& 2\chi'(t)\left( A\ppp_tu + \www{B}u + \int^t_0 (\ppp_tB)(t,\eta)ud\eta
\right)\\
-& \chi''(t)\left( Au + \int^t_0 B(t,\eta)u(x,\eta) d\eta\right)
\quad \mbox{in $Q$},\\
& \chi\ppp_t^2u\vert_{\ppp\Omega} = 0.
\end{array}\right.
                                              \eqno{(2.19)}
$$
Henceforth $U_2$ generically denotes a function satisfying
$$
U_2(x,t) \le C\sumalpha \sum_{k=0}^2 
\left( \vert \dalpha\ppp_t^ku\vert^2 
+ \left\vert \int^t_0 \vert \dalpha\ppp_t^ku\vert^2 d\eta\right\vert
\right)
$$
for $(x,t) \in Q$.  Applying Lemma 2.2 to (2.19) with $p=0$,
we obtain
\begin{align*}
&\gamma\int_Q \sumalpha \vert \dalpha(\chi\ppp_t^2u)\vert^2 \weight dxdt
\le C\gamma\int_Q s\va \vert \ppp_t^2v\vert^2 \weight dxdt\\
+& C\gamma\int_Q s\va\vert \chi \www{B}\ppp_tu\vert^2 \weight dxdt
+ C\gamma\int_Q s\va\chi^2 \vert (\ppp_t\www{B})(t)u 
+ (\ppp_tB)(t,t)u\vert^2 \weight dxdt \\
+& C\gamma\int_Q s\va\biggl(4 \vert \chi'\vert^2
\left\vert A\ppp_tu  + \www{B}u 
+ \int^t_0 (\ppp_tB)(t,\eta)ud\eta \right\vert^2\\
+ &\left\vert \chi''\left( Au + \int^t_0 Bud\eta \right)\right\vert^2
\biggr) \weight dxdt
+ Ce^{Cs}\Vert \ppp_{\nu}\ppp_t^2u\Vert^2_{L^2(\Sigma)}.
\end{align*}
We estimate the second term on the right-hand side as follows.
By $\ppp_tu(\cdot,0) = 0$, we have
$$
\chi\www{B}\ppp_tu(x,t) = \chi\www{B}\int^t_0 \ppp_t^2u(x,\eta) d\eta
= \chi\int^t_0 \www{B}\ppp_t^2u(x,\eta) d\eta.
$$
Lemma 2.3 with $q=1$ yields
\begin{align*}
& \int_Q s\gamma\va \vert \chi\www{B}\ppp_tu\vert^2 \weight dxdt
= \int_Q s\gamma\va\chi^2 \left\vert \int^t_0 \www{B}\ppp_t^2u(x,\eta) d\eta
\right\vert^2 \weight dxdt\\
\le& C\int_Q \chi^2 \vert \www{B}\ppp_t^2 u\vert^2 \weight dxdt
+ C\int_Q \vert \ppp_t(\chi^2)\vert \weight \left\vert 
\int^t_0 \vert \www{B}\ppp_t^2u(x,\eta)\vert^2 d\eta \right\vert dxdt\\
\le& C\int_Q \sumalpha \vert \dalpha(\chi\ppp_t^2u)\vert^2 \weight dxdt
+ C\int_Q \vert \ppp_t(\chi^2)\vert U_2(x,t)\weight dxdt.
\end{align*}
Consequently 
\begin{align*}
&\gamma\int_Q \sumalpha \vert \dalpha(\chi \ppp_t^2u)\vert^2 \weight dxdt\\
\le& C\int_Q s\gamma\va \vert \ppp_t^2v\vert^2 \weight dxdt
+ C\int_Q \sumalpha \vert \dalpha(\chi\ppp_t^2u)\vert^2 \weight dxdt\\
+& C\int_Q s\gamma\va \sumalpha \vert \dalpha(\chi u)\vert^2
\weight dxdt\\
+& C\int_Q s\gamma\va(\vert \ppp_t(\chi^2)\vert 
+ \vert \chi'(t)\vert^2 + \vert\chi''(t)\vert^2)
U_2\weight dxdt + Ce^{Cs}\Vert \ppp_{\nu}\ppp_t^2u\Vert^2_{L^2(\Sigma)}.
\end{align*}
Thus we can absorb the second term on the right-hand side into the
left-hand side, so that 
\begin{align*}
&\gamma\int_Q \sumalpha \vert \dalpha(\chi \ppp_t^2u)\vert^2 \weight dxdt\\
\le& C\int_Q s\gamma\va \vert \ppp_t^2v\vert^2 \weight dxdt
+ C\int_Q s\gamma\va \sumalpha \vert \dalpha(\chi u)\vert^2
\weight dxdt\\
+& C\int_Q s\gamma\va(\vert \ppp_t(\chi^2)\vert 
+ \vert \chi'(t)\vert^2 + \vert\chi''(t)\vert^2)
U_2\weight dxdt + Ce^{Cs}\Vert \ppp_{\nu}\ppp_t^2u\Vert^2_{L^2(\Sigma)}.
\end{align*}
Fixing $\gamma > 0$ sufficiently large and applying Theorem 1.1 (i) to 
the second term on the right-hand side, we obtain
$$
\int_Q \sumalpha \vert \dalpha(\chi \ppp_t^2u)\vert^2 \weight dxdt 
\le C\int_Q s\va \vert \ppp_t^2v\vert^2 \weight dxdt
                                               \eqno{(2.20)}
$$
\begin{align*}
+ &C\left(\int_Q \vert AF\vert^2 \weight dxdt
+ \Vert u\Vert^2_{H^1(-T,T;H^2(\OOO))} s^2\Phi^2e^{2s\delta}
+ e^{Cs}\Vert \ppp_{\nu}u\Vert^2_{H^2(-T,T;L^2(\Gamma))}\right)\\
+& C\int_Q s\va\mu(t) U_2\weight dxdt 
+ Ce^{Cs}\Vert \ppp_{\nu}\ppp_t^2u\Vert^2_{L^2(\Sigma)}\\
\le & C\int_Q s\va \vert \ppp_t^2v\vert^2 \weight dxdt
+ C \int_Q \vert AF\vert^2 \weight dxdt\\
+ & C\Vert u\Vert^2_{H^2(-T,T;H^2(\OOO))} s^2\Phi^2e^{2s\delta}
+ Ce^{Cs}\Vert \ppp_{\nu}u\Vert^2_{H^2(-T,T;L^2(\Gamma))}.
\end{align*}

Next, setting $v_1 = \ppp_tv$, we differentiate the first equation in (2.8)
in $t$, and we have
$$\left\{ \begin{array}{rl}
& \ppp_t^2v_1 - Av_1 = \chi A\ppp_tF + \chi'AF + \chi(\ppp_tJ)
+ \chi'J + \ppp_tS, \\
& \ppp_t^jv_1(\cdot,\pm T) = 0 \quad \mbox{in $\OOO$, $j=0,1$},\\
& v_1\vert_{\ppp\OOO} = 0.
\end{array}\right.           \eqno{(2.21)}
$$
Since
$$
\vert \chi\ppp_tJ(x,t)\vert^2 \le C\sumalpha \sum^2_{k=0}
\vert \dalpha(\chi \ppp_t^ku)\vert^2 \le CU_2(x,t), \quad 
(x,t) \in Q
                                         \eqno{(2.22)}
$$
and
$$
\vert \ppp_tS(x,t)\vert \le C\mu(t)U_2(x,t), \quad 
(x,t) \in Q,                                                      \eqno{(2.23)}
$$
in terms of (2.22) and (2.23), we apply Lemma 2.1 to (2.21) to 
obtain
\begin{align*}
& \int_Q (s\gamma\va\vert \nabla_{x,t}v_1\vert^2 + s^3\gamma^3\va^3
\vert v_1\vert^2) \weight dxdt
\le C\int_Q \vert \chi A\ppp_tF + \chi'AF\vert^2 \weight dxdt\\
+& C\int_Q \vert \chi\ppp_tJ + \chi'J + \ppp_tS\vert^2 \weight dxdt
+ Ce^{Cs}\Vert \ppp_{\nu}\ppp_tv\Vert^2_{L^2(\Sigma)}\\
\le& C\int_Q (\vert A\ppp_tF\vert^2 + \vert AF\vert^2) \weight dxdt
+ C\int_Q \sumalpha \sum_{k=0}^2 \vert\dalpha(\chi\ppp_t^ku)\vert^2
\weight dxdt\\
+& C\int_Q \mu(t)\sumalpha \sum_{k=0}^2 \vert \dalpha(\ppp_t^ku)\vert^2
\weight dxdt
+ Ce^{Cs}\Vert \ppp_{\nu}\ppp_tv\Vert^2_{L^2(\Sigma)}.
\end{align*}
Applying (2.20) and Theorem 1.1 (i) to the third and the second terms on the 
right-hand side respectively, we see
\begin{align*}
& \int_Q s\gamma\va\vert \ppp_t^2v\vert^2 \weight dxdt
\le C\int_Q (\vert A\ppp_tF\vert^2 + \vert AF\vert^2) \weight dxdt\\
+& C\Vert u\Vert^2_{H^2(-T,T;H^2(\OOO))}s^2\Phi^2e^{2s\delta}
+ Ce^{Cs}\Vert \ppp_{\nu}u\Vert^2_{H^2(-T,T;L^2(\Gamma))}\\
+& C\int_Q s\va\vert \ppp_t^2v\vert^2 \weight dxdt
+ Ce^{Cs}\Vert \ppp_{\nu}\ppp_tv\Vert^2_{L^2(\Sigma)}.
\end{align*}
By (2.4), we have $v=\chi\ppp_t^2u - \chi F$ and
\begin{align*}
& \ppp_{\nu}\ppp_tv = \chi'\ppp_{\nu}\ppp_t^2 u 
+ \chi\ppp_{\nu}\ppp_t^3u - \chi'\ppp_{\nu}F - \chi\ppp_{\nu}\ppp_tF\\
= &\chi'\ppp_{\nu}\ppp_t^2 u + \chi\ppp_{\nu}\ppp_t^3u
\end{align*}
on $\ppp\OOO$. Choosing $s>0$ and $\gamma>0$ sufficiently large,
we can absorb the fourth terms on the right-hand side into the left-hand side,
we have
$$
\int_Q (s\gamma\va\vert \ppp_t^2v_1\vert^2 + s^3\gamma^3\va^3
\vert v_1\vert^2) \weight dxdt                \eqno{(2.24)}
$$
$$
\le C\int_Q (\vert A\ppp_tF\vert^2 + \vert AF\vert^2) \weight dxdt
+ C\Vert u\Vert^2_{H^2(-T,T;H^2(\OOO))}s^2\Phi^2e^{2s\delta}
+ Ce^{Cs}\Vert \ppp_{\nu}u\Vert^2_{H^3(-T,T;L^2(\Gamma))}
$$
for $s > s_0$.

Substituting (2.24) into (2.20), we obtain
\begin{align*}
& \int_Q \sumalpha \vert \dalpha (\chi\ppp_t^2u)\vert^2
\weight dxdt
\le C\int_Q (\vert AF\vert^2 + \vert A\ppp_tF\vert^2) \weight 
dxdt\\
+& C\Vert u\Vert^2_{H^2(-T,T;H^2(\OOO))}s^2\Phi^2e^{2s\delta}
+ Ce^{Cs}\Vert \ppp_{\nu}u\Vert^2_{H^3(-T,T;L^2(\Gamma))}.
\end{align*}
Thus the proof of Theoerem 1.1 (ii) is completed.
\\

We close this section with the following lemma which is 
nothing but (2.24) where we fix $\gamma > 0$ sufficiently large.
The lemma plays an essential role for the proof of Theorem 1.3.
\\
\vspace{0.2cm}
{\bf Lemma 2.4.}\\
Let $u \in H^2(-T,T;H^2(\Omega))$ satify (1.16) and let (1.17) hold.
Under the same assumptions in Theorem 1.1 (ii), we have
\begin{align*}
& \int_Q (s\vert \nabla_{x,t}v_1\vert^2 + s^3\vert v_1\vert^2) 
\weight dxdt           
\le C\int_Q (\vert A((\ppp_tR)f)\vert^2 + \vert A(Rf)\vert^2) \weight dxdt\\
+ & C\Vert u\Vert^2_{H^2(-T,T;H^2(\OOO))}s^2\Phi^2e^{2s\delta}
+ Ce^{Cs}\Vert \ppp_{\nu}u\Vert^2_{H^3(-T,T;L^2(\Gamma))}
\end{align*}
for $s > s_0$.  Here we set
$$
v_1 = \ppp_tv, \quad 
v(x,t) = \chi(t)Au(x,t) 
+ \chi(t)\int^t_0 B(t,\eta)u(x,\eta) d\eta \quad \mbox{in $Q$}.
$$
\\
\vspace{0.2cm}

\section{Energy estimates}

For proving Theorems 1.2 and 1.3, we show energy estimates for hyperbolic
equations with integral terms.  Such an energy estimate is classical 
for hyperbolic equations without integral terms (e.g., Komornik
\cite{Ko}, Lions \cite{L}), but the presence of the integral terms
makes extra estimation demanded. 
\\
{\bf Lemma 3.1.}\\
We assume that 
$$\left\{\begin{array}{rl}
&a_{ij} = a_{ji} \in C(\ooo{\Omega}), \quad 1 \le i,j \le n, \\
&\mbox{there exists a constant $\mu_0 > 0$ such that}\\
&\sum_{i,j=1}^n a_{ij}(x)\xi_i\xi_j \ge \mu_0 \sum_{i=1}^n \xi_i^2,
\quad x \in \ooo{\Omega}, \thinspace \xi_1, ..., \xi_n \in \R.
\end{array}\right.                         \eqno{(3.1)}
$$
and
$$\left\{\begin{array}{rl}
&\nabla_{x,t}p_{\alpha} \in C(\ooo{\Omega}\times [0,T]^2), \quad
\vert \alpha\vert = 2, \\
&p_{\alpha} \in C(\ooo{\Omega}\times [0,T]^2), \quad
k=0,1,2, \thinspace \vert \alpha\vert \le 1. \\
\end{array}\right.              \eqno{(3.2)}
$$
If $u \in H^2(0,T;L^2(\Omega)) \cap H^1(0,T;H^1_0(\Omega)) \cap
L^2(0,T;H^2(\Omega))$ satisfies 
$$
\ppp_t^2u(x,t) = Au + \int^t_0 \sumalpha p_{\alpha}(x,t,\eta)
\dalpha u(x,\eta) d\eta + F(x,t), \quad x \in \Omega, \thinspace 
0 < t < T                     \eqno{(3.3)}
$$
and
$$
u\vert_{\ppp\Omega} = 0, \quad 0 < t < T,  \eqno{(3.4)}
$$
then there exists a constant $C>0$, which is independent of 
choice of $u$, such that 
$$
E(t) \le C(E(0) + \Vert F\Vert_{L^2(\Omega\times (0,T))}^2),
\quad 0 \le t \le T.              \eqno{(3.5)}
$$
Here and henceforth we set 
$$
E(u)(t) = E(t) 
= \int_{\Omega} \left( \vert \ppp_tu(x,t)\vert^2
+ \sumij a_{ij}(x)\ppp_iu(x,t)\ppp_ju(x,t)\right) dx, \quad 
0 \le t \le T.                         \eqno{(3.6)}
$$
\\
{\bf Proof of Lemma 3.1.}
We multiply (3.3) with 
$\ppp_t u(x,t)$ and integrate over $\Omega$:
By (3.1), (3.2) and (3.4) and integration by parts, we obtain
$$
\int_{\OOO} \ppp_t^2u(x,t)\ppp_tu(x,t) dx 
= \frac{1}{2}\ppp_t\int_{\OOO} \vert \ppp_tu(x,t)\vert^2 dx
$$
and
\begin{align*}
& -\int_{\OOO} \sumij \ppp_i(a_{ij}(x)\ppp_ju(x,t))\ppp_tu(x,t) dx
= \int_{\OOO} \sumij a_{ij}\ppp_ju(x,t) \ppp_i\ppp_tu(x,t) dx\\
=& \int_{\OOO} \left( \sum_{i>j} a_{ij}(x)
\ppp_t(\ppp_iu(x,t)\ppp_ju(x,t))
+ \frac{1}{2}\sum_{i=1}^n a_{ii}(x)\ppp_t(\vert \ppp_iu(x,t)\vert^2)
\right) dx\\
=&\frac{1}{2}\ppp_t\int_{\OOO}
\left( \sum_{i\ne j} a_{ij}\ppp_iu(x,t)\ppp_ju(x,t)
+ \sum_{i=1}^n a_{ii}\vert \ppp_iu(x,t)\vert^2\right) dx
= \frac{1}{2}\ppp_t\int_{\OOO}
\left( \sumij a_{ij}\ppp_iu(x,t)\ppp_ju(x,t) \right) dx.
\end{align*}
Next we calculate
$$
\int^t_0 p_{\alpha}(x,t,\eta)\dalpha u(x,\eta)\ppp_tu(x,t) d\eta.
$$
For $\vert \alpha\vert =2$, we set $\dalpha u = \ppp_i\ppp_ju$ with 
$1 \le i,j \le n$.
Noting 
$$
p_{ij}(x,t,\eta)\ppp_ju(x,\eta)\ppp_i\ppp_tu(x,t)
= \ppp_t(p_{ij}(x,t,\eta)\ppp_ju(x,\eta)\ppp_iu(x,t))
- \ppp_tp_{ij}(x,t,\eta)\ppp_ju(x,\eta)\ppp_iu(x,t),
$$ 
we write
\begin{align*}
& \int_{\OOO}\int^t_0 p_{ij}(x,t,\eta)\ppp_i\ppp_ju(x,\eta)\ppp_tu(x,t) 
d\eta dx\\
=& -\int_{\OOO}\int^t_0 \ppp_ip_{ij}(x,t,\eta)\ppp_ju(x,\eta)\ppp_tu(x,t)
d\eta dx
- \int_{\OOO}\int^t_0 p_{ij}(x,t,\eta)\ppp_ju(x,\eta)\ppp_t\ppp_iu(x,t) 
d\eta dx\\
=& \biggl\{-\int_{\OOO}\int^t_0 \ppp_ip_{ij}(x,t,\eta)\ppp_ju(x,\eta)
\ppp_tu(x,t) d\eta dx
+ \int_{\OOO}\int^t_0 (\ppp_tp_{ij})(x,t,\eta)\ppp_ju(x,\eta)\ppp_iu(x,t) 
d\eta dx\biggr\}\\
- &\int_{\OOO}\int^t_0 \ppp_t(p_{ij}(x,t,\eta)\ppp_ju(x,\eta)\ppp_iu(x,t))
d\eta dx\\
=:& \www{J} - \int_{\OOO}\int^t_0 \ppp_t(p_{ij}(x,t,\eta)\ppp_ju(x,\eta)
\ppp_iu(x,t))d\eta dx.
\end{align*}
Here
\begin{align*}
& -\int_{\OOO}\int^t_0 \ppp_t(p_{ij}(x,t,\eta)\ppp_ju(x,\eta)\ppp_iu(x,t))
d\eta dx\\
=& - \int_{\OOO} \left( \ppp_t\left(\int^t_0 p_{ij}(x,t,\eta)\ppp_ju(x,\eta)
\ppp_iu(x,t)d\eta\right)
- p_{ij}(x,t,t)\ppp_ju(x,t)\ppp_iu(x,t)\right) dx.
\end{align*}
Therefore
$$
\int_{\OOO}\int^t_0 p_{ij}(x,t,\eta)\ppp_i\ppp_ju(x,\eta)\ppp_iu(x,t)
d\eta dx
=: J_{ij} - \ppp_t\int^t_0\int_{\OOO}p_{ij}(x,t,\eta)\ppp_ju(x,\eta)
\ppp_iu(x,t) dxd\eta,                   \eqno{(3.7)}
$$
where
$$
J_{ij} = -\int_{\OOO}\int^t_0 \ppp_ip_{ij}(x,t,\eta)\ppp_ju(x,\eta)
\ppp_tu(x,t) d\eta dx                       \eqno{(3.8)}
$$
$$
+ \int_{\OOO}\int^t_0 \ppp_tp_{ij}(x,t,\eta)\ppp_ju(x,\eta)\ppp_iu(x,t)
d\eta dx
+ \int_{\OOO} p_{ij}(x,t,t)\ppp_ju(x,t)\ppp_iu(x,t) dx.
$$
For $\vert \alpha \vert =1$, we have
$$
\left\vert \int_{\OOO}\int^t_0 p_i(x,t,\eta)\ppp_ju(x,\eta)\ppp_tu(x,t) dx
\right\vert                                  \eqno{(3.9)}
$$
\begin{align*}
& \le C\int^t_0 \left(\int_{\OOO}
(\vert \ppp_ju(x,t,\eta)\vert^2 + \vert \ppp_tu(x,t)\vert^2) dx
\right) d\eta\\
\le &C\int^t_0 (E(\eta) + E(t)) d\eta \le CE(t) + C\int^t_0 E(\eta) d\eta
\end{align*}
and
for $\vert \alpha \vert = 0$, we use the Poincar\'e inequality to obtain
$$
\left\vert \int_{\OOO}\int^t_0 p_0(x,t,\eta)u(x,\eta)\ppp_iu(x,t)d\eta dx
\right\vert                           
                                   \eqno{(3.10)}
$$
\begin{align*}
\le &C\int^t_0\int_{\OOO} (\vert \ppp_tu(x,t)\vert^2
+ \vert u(x,\eta)\vert^2) dx d\eta
\le C\int^t_0 (E(\eta) + E(t)) d\eta \\
\le & CE(t) + C\int^t_0 E(\eta) d\eta.
\end{align*}
Hence
\begin{align*}
& \int_{\OOO}\int^t_0 \sumalpha p_{\alpha}(x,t,\eta)\dalpha u(x,\eta)
\ppp_tu(x,t)) d\eta dx\\
=& \sumij J_{ij} 
- \ppp_t \int^t_0 \sumij \int_{\OOO} p_{ij}(x,t,\eta)\ppp_ju(x,\eta)
\ppp_iu(x,t) dxd\eta \\
+ &\sum_{j=1}^n \int_{\OOO}\int^t_0 p_j(x,t,\eta)\ppp_ju(x,\eta)\ppp_iu(x,t)
d\eta dx\\
+& \int_{\OOO}\int^t_0 p_0(x,t,\eta)u(x,\eta)\ppp_tu(x,t) d\eta dx.
\end{align*}
Thus we obtain
\begin{align*}
& \frac{1}{2}\frac{dE}{dt}(t) 
=  \sumij J_{ij} 
- \ppp_t \int^t_0 \sumij \int_{\OOO} p_{ij}(x,t,\eta)\ppp_ju(x,\eta)
\ppp_iu(x,t) dxd\eta \\
+ & \sum_{j=1}^n \int_{\OOO}\int^t_0 p_j(x,t,\eta)\ppp_ju(x,\eta)\ppp_iu(x,t)
d\eta dx
+ \int_{\OOO}\int^t_0 p_0(x,t,\eta)u(x,\eta)\ppp_tu(x,t) d\eta dx.
\end{align*}
We choose $r \in (0,T)$ arbitrarily.
Integrating both sides with respect to $t$ in $(0,r)$,
we reach
\begin{align*}
& \frac{1}{2}(E(r) - E(0)) 
= \sumij \int^r_0 J_{ij}(t) dt 
- \int^r_0 \sumij \int_{\OOO} p_{ij}(x,r,\eta)\ppp_ju(x,\eta)
\ppp_iu(x,t) dxd\eta \\
+ & \sum_{j=1}^n \int^r_0 \int_{\OOO}\int^t_0 
p_j(x,t,\eta)\ppp_ju(x,\eta)\ppp_iu(x,t) d\eta dx dt
+ \int^r_0 \int_{\OOO}\int^t_0 p_0(x,t,\eta)u(x,\eta)\ppp_tu(x,t) 
d\eta dx dt\\
+& \int^r_0\int_{\OOO} F(x,t)\ppp_tu(x,t) dxdt.
\end{align*}
Here we have
\begin{align*}
& \left\vert \sumij \int^r_0 J_{ij}(t) dt \right\vert 
\le \sumij \int^r_0 \int_{\OOO} \int^t_0
\vert \ppp_ip_{ij}(x,t,\eta)\vert \vert \ppp_ju(x,\eta)\vert
\vert \ppp_tu(x,t)\vert d\eta dx dt \\ 
+& \sumij \int^r_0 \int_{\OOO} \int^t_0
\vert \ppp_tp_{ij}(x,t,\eta)\vert \vert \ppp_ju(x,\eta)\vert
\vert \ppp_iu(x,t)\vert d\eta dx dt \\
+ & \sumij \int^r_0 \int_{\OOO} \vert p_{ij}(x,t,t)\vert 
\vert \ppp_ju(x,t)\vert \vert \ppp_iu(x,t)\vert dx dt \\
\le& C\int^r_0 \int^t_0 (E(\eta) + E(t)) d\eta dt 
+ C\int^r_0 E(t) dt\\
\le& C\int^r_0 \left( \int^r_{\eta} dt \right) E(\eta) d\eta 
+ C\int^r_0 tE(t) dt + C\int^t_0 E(t) dt
\le C\int^r_0 E(t) dt.
\end{align*}
Moreover, for any $\ep > 0$, we can choose a constant $C_{\ep} > 0$ 
such that 
\begin{align*}
& \left\vert  \sumij \int^r_0 \int_{\OOO} 
\ppp_ip_{ij}(x,r,\eta) \ppp_ju(x,\eta) \ppp_iu(x,r) dx d\eta \right\vert
\le C\int^r_0 \sumij \int_{\OOO} 
\vert \ppp_ju(x,\eta)\vert \vert \ppp_iu(x,r)\vert dx d\eta\\
\le& C\int^r_0 \int_{\OOO} \vert \nabla u(x,\eta)\vert 
\vert \nabla u(x,r) \vert dxd\eta 
\le C_1\ep\int^r_0 \int_{\OOO} \vert \nabla u(x,r)\vert^2 dxd\eta
+ C_{\ep}\int^r_0 \int_{\OOO} \vert \nabla u(x,\eta)\vert^2 dxd\eta\\
\le& C_1\ep T E(r) + C_{\ep}\int^r_0 E(\eta) d\eta.
\end{align*}
Moreover (3.9) and (3.10) yield
\begin{align*}
&\biggl\vert \sum_{j=1}^n \int^r_0 \int_{\OOO} \int^t_0
p_j(x,t,\eta) \ppp_ju(x,\eta) \ppp_tu(x,t) d\eta dxdt\\
+ &\int^r_0\int_{\OOO} \int^t_0 p_0(x,t,\eta) u(x,\eta)\ppp_tu(x,t) 
d\eta dxdt \biggr\vert
\le C\int^r_0 E(t) dt
\end{align*}
and
\begin{align*}
&\left\vert \int^r_0 \int_{\OOO} F(x,t)\ppp_tu(x,t) dxdt \right\vert
\le \int^r_0 \int_{\OOO} \vert F(x,t)\vert^2 dxdt 
+ \int^r_0 \int_{\OOO} \vert \ppp_tu(x,t)\vert^2 dxdt\\
\le& \Vert F\Vert^2_{L^2(\Omega \times (0,T))}
+ C\int^r_0 E(t) dt.
\end{align*}
Thus we obtain
$$
E(r) - E(0) \le C_2\ep E(r) + C_{\ep}\int^r_0 E(t) dt 
+ C\Vert F\Vert^2_{L^2(\OOO\times (0,T))}.     \eqno{(3.11)}
$$
Choosing $\ep > 0$ sufficiently small such that $C_2\ep < \frac{1}{2}$,
we have
$$
\frac{1}{2}E(r) \le E(0) + C\int^r_0 E(t) dt + C\Vert F\Vert^2
_{L^2(\OOO\times (0,T))}.
$$
The Gronwall inequality implies 
$$
E(r) \le C(E(0) + \Vert F\Vert^2_{L^2(\OOO\times (0,T))})e^{CT},
\quad 0 \le r \le T.
$$
Thus the proof of Lemma 3.1 is completed.
\\
\vspace{0.3cm}

Next, on the basis of Lemma 3.1,  
we prove an energy estimate in Sobolev spaces of higher orders for
the solution to (1.16).
We assume that $u \in H^2(0,T;H^2(\OOO))$ satisfies (1.16),
$\ppp_{\nu}u \in H^3(0,T;L^2(\Gamma))$ and that (1.17) holds.
Then we prove
\\
{\bf Lemma 3.2.}\\
There exists a constant $C>0$ such that 
$$
\sum_{k=0}^2 \Vert \ppp_t^ku(\cdot,t)\Vert_{H^2(\OOO)}
\le C\Vert f\Vert_{H^2(\OOO)}, \quad 0 \le t \le T
$$
for each $f \in H^2(\OOO) \cap H^1_0(\OOO)$.
\\
\vspace{0.2cm}
{\bf Proof of Lemma 3.2.}
We set 
$$
u_1 = \ppp_tu, \quad u_2 = \ppp_t^2u, \quad u_3 = \ppp_t^3u
                                                     \eqno{(3.12)}
$$
and
$$\left\{ \begin{array}{rl}
& b_{\alpha}^{(1)}(x,t,\eta) = b_{\alpha}(x,t,t)
+ \int^t_{\eta} \ppp_tb_{\alpha}(x,t,\xi) d\xi,\\
& b_{\alpha}^{(k+1)}(x,t,\eta) = b_{\alpha}^{(k)}(x,t,t)
+ \int^t_{\eta} \ppp_tb_{\alpha}^{(k)}(x,t,\xi) d\xi,\quad k=1,2.
\end{array}\right.                           \eqno{(3.13)}
$$
We set 
$$
c_{\alpha}(x,t) = b_{\alpha}^{(3)}(x,t,0).
$$
Then 
$$\left\{ \begin{array}{rl}
& \ppp_t^2u_1(x,t) = Au_1 + \int^t_0 \sumalpha b_{\alpha}^{(1)}(x,t,\eta)
\dalpha u_1(x,\eta) d\eta + (\ppp_tR)f, \\
&u_1(x,0) = 0, \quad \ppp_tu_1(x,0) = R(x,0)f(x), \quad x\in \OOO,\\
& u_1\vert_{\ppp\OOO} = 0,
\end{array}\right.                        \eqno{(3.14)}
$$
$$\left\{ \begin{array}{rl}
& \ppp_t^2u_2(x,t) = Au_2 + \int^t_0 \sumalpha b_{\alpha}^{(2)}(x,t,\eta)
\dalpha u_2(x,\eta) d\eta + (\ppp_t^2R)f, \\
&u_2(x,0) = R(x,0)f(x), \quad \ppp_tu_2(x,0) = \ppp_tR(x,0)f(x), 
\quad x\in \OOO,\\
& u_2\vert_{\ppp\OOO} = 0
\end{array}\right.                        \eqno{(3.15)}
$$
and
$$\left\{ \begin{array}{rl}
& \ppp_t^2u_3(x,t) = Au_3 + \int^t_0 \sumalpha b_{\alpha}^{(3)}(x,t,\eta)
\dalpha u_3(x,\eta) d\eta + (\ppp_t^3R)f \\
+ &\sumalpha c_{\alpha}(x,t)\dalpha (R(x,0)f), \\
&u_3(x,0) = (\ppp_tR(x,0))f, \quad 
\ppp_tu_3(x,0) = A(R(x,0)f), \quad x\in \OOO,\\
& u_3\vert_{\ppp\OOO} = 0.
\end{array}\right.                        \eqno{(3.16)}
$$
Indeed we can verify (3.14) as follows.  We differentiate (1.16) in $t$ to have
$$
\ppp_t^2u_1(x,t) = Au_1 + \sumalpha b_{\alpha}(x,t,t)\dalpha u(x,t)
+ \int^t_0 \sumalpha \ppp_tb_{\alpha}(x,t,\eta)\dalpha u(x,\eta) d\eta 
+ (\ppp_tR)f.
$$
Noting that $u(x,t) = \int^t_0 u_1(x,\xi) d\xi$ by $u(\cdot,0) = 0$ and 
changing the orders of integration, we obtain
$$\left\{ \begin{array}{rl}
& \ppp_t^2u_1(x,t) = Au_1 + \sumalpha \int^t_0 b_{\alpha}(x,t,t)
\dalpha u_1(x,\eta) d\eta \\
+ &\sumalpha \int^t_0 (\ppp_tb_{\alpha})(x,t,\eta)
\left(\int^{\eta}_0 \dalpha u_1(x,\xi) d\xi \right) d\eta
+ (\ppp_tR)f \\
= & Au_1 + \int^t_0 \sumalpha b_{\alpha}(x,t,t)
\dalpha u_1(x,\eta) d\eta 
+ \int^t_0 \left( \int^t_{\xi} (\ppp_tb_{\alpha})(x,t,\eta) d\eta \right)
\dalpha u_1(x,\xi) d\xi + (\ppp_tR)f 
\end{array}\right.                        
$$
and
$$
u_1(x,0) = \ppp_tu(x,0) = 0,
$$
\begin{align*}
& \ppp_tu_1(x,0) = \ppp_t^2u(x,0) \\
= &\left( Au + \int^t_0 \sumalpha b_{\alpha}(x,t,\eta)
\dalpha u(x,\eta) d\eta\right)\biggl\vert_{t=0}
+ R(x,0)f(x) = R(x,0)f.
\end{align*}
Therefore (3.14) is verified.  The systems (3.15) and (3.16) can be  
verified similarly by noting that 
$$
u_1(x,t) = \int^t_0 u_2(x,\xi) d\xi, \quad
u_2(x,t) = \int^t_0 u_3(x,\xi) d\xi + R(x,0)f(x).
$$

Applying Lemma 3.1 to (3.14) - (3.16), we have
$$
\sum_{k=0}^3 (\Vert \ppp_t^{k+1} u(\cdot,t)\Vert_{L^2(\OOO)}
+ \Vert \ppp_t^ku(\cdot,t)\Vert_{H^1(\OOO)})
\le C\Vert f\Vert_{H^2(\OOO)}, \quad 0 \le t \le T.   
                                                  \eqno{(3.17)}
$$
Next we have to estimate $\Vert \ppp_t^j u(\cdot,t)\Vert_{H^2(\OOO)}$,
$k=0,1,2$.  Since $\ppp_tu(\cdot,t) = \int^t_0 \ppp_t^2u(\cdot,\eta)d\eta$
and $u(\cdot,t) = \int^t_0 (t-\eta)\ppp_t^2u(\cdot,\eta)d\eta$ by
$u(\cdot,0) = \ppp_tu(\cdot,0) = 0$ in $\OOO$, it suffices to estimate
$\Vert \ppp_t^2u(\cdot,t)\Vert_{H^2(\OOO)}$.
By (3.17) we have
$$
\Vert \ppp_t^4u(\cdot,t)\Vert_{L^2(\OOO)} \le C\Vert f\Vert_{H^2(\OOO)},
\quad 0 \le t \le T.                   \eqno{(3.18)}
$$
Therefore (3.15) implies
$$
A\ppp_t^2u = \ppp_t^4u(x,t) - \int^t_0 \sumalpha b_{\alpha}^{(2)}
(x,t,\eta) \dalpha \ppp_t^2u(x,\eta) d\eta
+ (\ppp_t^2R)f,
$$
and so
$$
\Vert A\ppp_t^2u(\cdot,t)\Vert_{L^2(\OOO)}
\le C\Vert \ppp_t^4u(\cdot,t)\Vert_{L^2(\OOO)}
+ C\int^t_0 \Vert \ppp_t^2u(\cdot,\eta)\Vert_{H^2(\OOO)} d\eta
+ C\Vert f\Vert_{H^2(\OOO)}.
$$
Since $\ppp_t^2u(\cdot,t)\vert_{\ppp\OOO} = 0$, we apply the a priori
estimate for the elliptic boundary value problem, by (3.18) we obtain
$$
\Vert \ppp_t^2u(\cdot,t)\Vert_{H^2(\OOO)}
\le C\Vert f\Vert_{H^2(\OOO)} + C\int^t_0 \Vert \ppp_t^2u(\cdot,\eta)\Vert
_{H^2(\OOO)} d\eta, \quad 0 \le t \le T.
$$
The Gronwall inequality yields
$$
\Vert \ppp_t^2u(\cdot,t)\Vert_{H^2(\OOO)}
\le C\Vert f\Vert_{H^2(\OOO)}, \quad 0 \le t \le T.
$$
Thus the proof of Lemma 3.2 is completed.

In this section we further show the following lemma for the hyperbolic
equation without integral terms.
\\
{\bf Lemma 3.3.}\\
We assume (3.1).  Let $w \in H^2(0,T;L^2(\OOO)) \cap L^2(0,T;H^2(\OOO))
\cap C([0,T];H^1_0(\OOO)) \cap C^1([0,T];L^2(\OOO))$ satisfy
$$
\ppp_t^2 w(x,t) = Aw(x,t) + G(x,t), \quad x \in \OOO, \thinspace
0 < t < T.                                     \eqno{(3.19)}
$$
\\
(i) There exists a constant $C>0$ such that 
$$
\Vert \ppp_{\nu}w\Vert^2_{L^2(0,T;L^2(\ppp\OOO))}
\le C(E(0) + \Vert G\Vert^2_{L^2(0,T;L^2(\OOO))}).   \eqno{(3.20)}
$$
(ii) There exists a constant $C >0$ such that 
$$
E(0) \le C\left(E(t) + \int^t_0 \Vert G(\cdot,\eta)\Vert^2
_{L^2(\OOO)} d\eta \right), \quad 0 \le t \le T.
                                                    \eqno{(3.21)}
$$
\\
{\bf Proof of Lemma 3.3.} The estimate (3.20) is proved by the multiplier 
method (e.g., Komornik \cite{Ko}).  
That is, we choose $h: \OOO \longrightarrow {\R}^n$ 
such that $h \in C^1(\overline{\OOO})$ and $h\vert_{\ppp\OOO} = \nu$ which is 
the unit outward normal vector to $\ppp\OOO$.  Then multiplyng (3.19) with 
$h\cdot\nabla w$ and integrating over $\OOO\times (0,T)$, we see 
(3.20).  We omit the details and see e.g., \cite{Ko} for the complete proof.

A usual energry estimate yields
$$
E(w)(t) \le C(E(w)(0 ) + \Vert G\Vert^2_{L^2(0,t;L^2(\Omega))}).  
                                                 \eqno{(3.22)}
$$
Since (3.19) is time-reversing, we can consider (3.19) by regarding 
$t$ as initial time and, applying (3.22), we obtain (3.21).
Thus the proof of Lemma 3.3 is completed.
\section{Proof of Theorem 1.2}

Our proof is a modification of Kazemi and Klibanov \cite{KK} and 
Klibanov and Malinsky \cite{KM} which discuss hyperbolic equations 
without integral terms.

We make the even extension of $y$ to $(-T,0)$:
$$
y(x,t) = 
\left\{ \begin{array}{rl}
y(x,t), \quad & 0 < t < T, \\
y(x,-t), \quad & -T<t<0.
\end{array}\right.
                                \eqno{(4.1)}
$$
Then, by $\ppp_ty(\cdot,0) = 0$, we can verify that
$$
\ppp_ty(x,t) = 
\left\{ \begin{array}{rl}
\ppp_ty(x,t), \quad & 0 < t < T, \\
-\ppp_ty(x,-t), \quad & -T<t<0, 
\end{array}\right.
\quad 
\ppp_t^2y(x,t) = 
\left\{ \begin{array}{rl}
\ppp_t^2y(x,t), \quad & 0 < t < T, \\
\ppp_t^2y(x,-t), \quad & -T<t<0,
\end{array}\right.      
$$
$$
\ppp_t^3y(x,t) = 
\left\{ \begin{array}{rl}
\ppp_t^3y(x,t), \quad & 0 < t < T, \\
-\ppp_t^3y(x,-t), \quad & -T<t<0.
\end{array}\right.
                                      \eqno{(4.2)}
$$
Hence
$$\left\{ \begin{array}{rl}
& y \in C^2([-T,T];H^2(\OOO)\cap H^1_0(\OOO)), \\
& \dalpha y \in H^2(-T,T;L^2(\OOO)) \cap L^2(-T,T;H^2(\OOO)).
\end{array}\right.
                                                \eqno{(4.3)}
$$
Next we will estimate $\Vert y\Vert_{H^1(-T,T;H^2(\OOO))}$.
By (4.1) and (4.2), it is sufficient to estimate for $0 < t < T$, that is,
$\Vert y\Vert_{H^1(0,T;H^2(\OOO))}$.
Similarly to the proof of Lemma 3.2, we set 
$$
y_1 = \ppp_ty, \quad y_2 = \ppp_t^2y,
$$
\begin{align*}
&b_{\alpha}^{(1)}(x,t,\eta) = b_{\alpha}(x,t,t)
+ \int^t_{\eta} \ppp_tb_{\alpha}(x,t,\xi)d\xi,\\
&b_{\alpha}^{(2)}(x,t,\eta) = b_{\alpha}^{(1)}(x,t,t)
+ \int^t_{\eta} \ppp_tb_{\alpha}^{(1)}(x,t,\xi)d\xi.
\end{align*}
We recall that 
$$
a(x) := y(x,0), \quad x \in \Omega.
$$
Noting that $y(x,t) = \int^t_0 y_1(x,\xi) d\xi + a(x)$ and 
$y_1(x) = \int^t_0 y_2(x,\xi)d\xi$ and 
$y_1(x,0) = \ppp_ty(x,0) =0$, in a way similar to (3.13) and (3.14), we can 
$$\left\{ \begin{array}{rl}
& \ppp_t^2y_1(x,t) = Ay_1 + \int^t_0 \sumalpha b_{\alpha}^{(1)}(x,t,\eta)
\dalpha y_1(x,\eta) d\eta + \sumalpha b_{\alpha}^{(1)}(x,t,0)\dalpha a(x), \\
&y_1(x,0) = 0, \quad \ppp_ty_1(x,0) = Aa(x), \quad x\in \OOO,\\
& y_1\vert_{\ppp\OOO} = 0
\end{array}\right.                        \eqno{(4.4)}
$$
and
$$\left\{ \begin{array}{rl}
& \ppp_t^2y_2(x,t) = Ay_2 + \int^t_0 \sumalpha b_{\alpha}^{(2)}(x,t,\eta)
\dalpha y_2(x,\eta) d\eta + \sumalpha \ppp_tb_{\alpha}^{(1)}(x,t,0)
\dalpha a(x), \\
&y_2(x,0) = Aa(x), \quad \ppp_ty_2(x,0) = \sumalpha b_{\alpha}^{(1)}(x,0,0)
\dalpha a(x), \quad x\in \OOO,\\
& y_2\vert_{\ppp\OOO} = 0.          
\end{array}\right.             \eqno{(4.5)}
$$
Applying Lemma 3.1 to (4.4) and (4.5), we obtain
$$
\Vert \ppp_t^3y(\cdot,t)\Vert_{L^2(\OOO)} 
\le \sum_{k=0}^1 (\Vert \ppp_ty_k(\cdot,t)\Vert_{L^2(\OOO)}
+ \Vert \nabla y_k(\cdot,t)\Vert_{L^2(\OOO)})
\le C\Vert a\Vert_{H^2(\OOO)}, \quad 0 \le t \le T.               \eqno{(4.6)}
$$
By the first equation in (4.4), we have
$$
A\ppp_ty(\cdot,t) = \ppp_t^3y - \sumalpha \int^t_0 b_{\alpha}^{(1)}(x,t,\eta)
\dalpha \ppp_ty(x,\eta) d\eta - \sumalpha b_{\alpha}^{(1)}(x,t,0)\dalpha a(x)
$$
and
$$
\ppp_ty(\cdot,t)\vert_{\ppp\OOO} = 0.
$$
Applying the a priori estimate for the elliptic boundary value problem,
we reach 
$$
\Vert \ppp_ty(\cdot,t)\Vert_{H^2(\OOO)}
\le C\Vert \ppp_t^3y(\cdot,t)\Vert_{L^2(\OOO)}
+ C\int^t_0 \Vert \ppp_ty(x,\eta)\Vert_{H^2(\OOO)} d\eta 
+ C\Vert a\Vert_{H^2(\OOO)}.
$$
>From (4.6) and the Gronwall inequality it follows that
$$
\Vert \ppp_ty(\cdot,t)\Vert_{H^2(\OOO)}
\le C\Vert a\Vert_{H^2(\OOO)}, \quad 0 \le t \le T.             \eqno{(4.7)}
$$
Since $y(\cdot,t) = \int^t_0 \ppp_ty(\cdot,\xi) d\xi + a$, it follows 
from (4.7) that 
$$
\Vert y(\cdot,t)\Vert_{H^2(\OOO)} \le C\Vert a\Vert_{H^2(\OOO)},
\quad 0 \le t \le T.
$$
Thus 
$$
\Vert y(\cdot,t)\Vert_{H^2(\OOO)} + \Vert \ppp_ty(\cdot,t)\Vert
_{H^2(\OOO)} \le C\Vert a\Vert_{H^2(\OOO)}, \quad 0 \le t \le T.
                                                       \eqno{(4.8)}
$$

Now we choose $\ep > 0$ and $\delta > 0$ in (1.10). 
The assumption (1.15) implies 
$$
\psi(x,\pm T) = \vert x-x_0\vert^2 - \beta T^2 < 0, \quad 
x\in \ooo{\OOO},                                            \eqno{(4.9)}
$$
and
$$
\psi(x,0) = \vert x-x_0\vert > 0, \quad x\in \ooo{\OOO}   
$$
by $x_0 \not\in \ooo{\OOO}$.  Therefore there exist $\ep_0 > 0$ and
$\ep \in \left(0, \frac{T}{4}\right)$ such that
$$\left\{ \begin{array}{rl}
\va(x,t) \le 1 - \ep_0, \quad 
&x\in \ooo{\OOO}, \thinspace \vert T-t\vert \le 2\ep \thinspace
\mbox{or} \thinspace \vert T+t\vert \le 2\ep,\\
\va(x,t) \ge 1+ \ep_0, \quad &x \in \ooo{\OOO}, \thinspace
\vert t\vert \le 2\ep.
\end{array}\right.                       \eqno{(4.10)}
$$
By (4.10) we note that $\va(x,\pm T) < 1$ for $x \in \overline{\OOO}$, and 
so $0 < \delta < 1$ if we choose $\ep>0$ sufficiently small in (1.10).
Then we can set $\delta = 1 - \ep_0$ with $\ep_0 > 0$.  Then
$$
\va(x,t) \ge 1 + \ep_0 = \delta + 2\ep_0, \quad \vert t\vert \le 2\ep.
                                                            \eqno{(4.11)}
$$

We set 
$$
z(x,t) = \chi(t)\ppp_t^2y(x,t), \quad (x,t) \in Q.
$$
Noting that $F=0$ in $Q$, we apply (1.11) in Theorem 1.1 to 
(1.14) where we fix $\gamma > 0$ large and applying (4.8) and so
$$
\int^T_0 \int_{\OOO} (\vert \ppp_tz\vert^2 + \vert \nabla z\vert^2)
\weight dxdt                        \eqno{(4.12)}
$$
$$
\le C\Vert a\Vert^2_{H^2(\OOO)}s^2e^{2s\delta}
+ Ce^{Cs}\Vert \ppp_{\nu}y\Vert^2_{H^2(-T,T;L^2(\Gamma))}
$$
for all large $s>0$.  On the other hand, similarly to (2.8), we see
$$\left\{ \begin{array}{rl}
& \ppp_t^2 z(x,t) = Az + \widetilde{S}(x,t), \quad (x,t) \in Q, \\
& z(x,0) = (Aa)(x), \quad \ppp_tz(x,0) = \sumalpha b_{\alpha}(x,0,0)
\dalpha a(x), \\
& z\vert_{\ppp\OOO\times (-T,T)} = 0,
\end{array}\right.
                              \eqno{(4.13)}
$$
where 
$$
\vert \www{S}(x,t)\vert \le C\left( \sumalpha \sum_{k=0}^1
\vert \dalpha \ppp_t^ky(x,t)\vert
+ \int^t_0 \sumalpha \vert \dalpha y(x,\eta)\vert d\eta\right)
                                                  \eqno{(4.14)}
$$
for $(x,t) \in Q$.  Here, in terms of $\chi(0) = 1$, 
$\chi'(0) = 0$ and the first equation in (1.14), 
we calculated $z(x,0) = \ppp_t^2y(x,0)$ and 
$\ppp_tz(x,0) = \ppp_t^3y(x,0)$.

Applying (3.21) in Lemma 3.3, we fix sufficiently small $\delta_0 > 0$ 
such that $0 < \delta_0 < 2\ep$ and 
$$
E(0) \le C\left( E(t) + \int^t_0 \Vert \www{S}(\cdot,\eta)
\Vert^2_{L^2(\OOO)} d\eta \right), \quad 0 \le t \le \delta_0.
                                                     \eqno{(4.15)}
$$
Here we recall that $E(z)(t) = E(t) = \int_{\OOO}
\left( \vert \ppp_tz(x,t)\vert^2 + \sumij a_{ij}(x)
\ppp_iz(x,t)\ppp_jz(x,t)\right) dx$.  By (4.14) and (4.8), we obtain
$$
\int^t_0 \Vert \www{S}(\cdot,\eta)\Vert^2_{L^2(\OOO)} d\eta 
                                           \eqno{(4.16)}
$$
\begin{align*}
\le& C\int^t_0 \left( \sum_{k=0}^1 \Vert \ppp_t^ky(\cdot,\eta)\Vert
^2_{H^2(\OOO)} 
+ \int^{\eta}_0 \Vert y(\cdot,\xi)\Vert
^2_{H^2(\OOO)} d\xi\right) d\eta\\
\le& C\int^t_0 \sum_{k=0}^1 \Vert \ppp_t^ky(\cdot,\eta)\Vert
^2_{H^2(\OOO)} d\eta \le Ct\Vert a\Vert^2_{H^2(\OOO)}
\le CtE(0), \quad 0 \le t \le T.
\end{align*} 
At the last inequality, in view of the Poincar\'e inequality,  
$z(\cdot,0) = Aa$ and the a priori estimate for the boundary 
value problem for $A$, we used 
$$
E(0) \ge \int_{\OOO} \sumij a_{ij}\ppp_iz(x,0)\ppp_jz(x,0) dx
\ge C\int_{\OOO} \vert \nabla z(x,0)\vert^2 dx
                                                \eqno{(4.17)}
$$
$$
\ge C\int_{\OOO} \vert z(x,0)\vert^2 dx = C\int_{\OOO} \vert (Aa)(x)\vert^2 dx
\ge C\Vert a\Vert^2_{H^2(\OOO)}.
$$

Hence, by (4.15), we have 
$$
(1-C\delta_0)E(0) \le CE(t), \quad 0\le t\le \delta_0.
$$
We further choose small $\delta_0 > 0$ such that $1-C\delta_0 > 0$ and 
fix.  Then $E(0) \le C_1E(t)$ for $0 \le t \le \delta_0$.
Therefore
\begin{align*}
& C_1E(0)\delta_0e^{2s(\delta+2\ep_0)}
\le C_1\int^{\delta_0}_0 E(t)\weight dxdt\\
\le & \int^{\max\{2\ep,\delta_0\}}_0\int_{\OOO}
(\vert \ppp_tz\vert^2 + \vert \nabla z\vert^2)\weight dxdt
\le \int_Q (\vert \ppp_tz\vert^2 + \vert \nabla z\vert^2)\weight dxdt
\end{align*}
for all large $s>0$.
Substituting this into the left-hand side of (4.12) and using (4.17),
we obtain
$$
C_1E(0)\delta_0 e^{2s(\delta+2\ep_0)}
\le CE(0)s^2e^{2s\delta} + Ce^{Cs}\Vert \ppp_{\nu}y\Vert^2
_{H^2(-T,T;L^2(\Gamma))}
$$
for all large $s>0$.  Hence
$$
(C_1\delta_0 - Cs^2e^{-4s\ep_0})E(0) 
\le Ce^{Cs}\Vert \ppp_{\nu}y\Vert^2_{H^2(-T,T;L^2(\Gamma))}
$$
for all large $s>0$.  Choosing $s>0$ sufficiently large such that 
$C_1\delta_0 - Cs^2e^{-4s\ep_0}>0$, we complete the proof of the second
inequality in the conclusion of Theorem 1.2.

Next we prove the first inequality of the conclusion.
In place of $z = \chi\ppp_t^2y$, we set $y_2 = \ppp_t^2y$.
Then, similarly to (2.8), we have
$$\left\{ \begin{array}{rl}
& \ppp_t^2 y_2(x,t) = Ay_2 + \widetilde{S_1}(x,t), \quad (x,t), \quad
x \in \OOO, \thinspace 0 < t < T, \\
& y_2(x,0) = (Aa)(x), \quad \ppp_ty_2(x,0) = \sumalpha b_{\alpha}(x,0,0)
\dalpha a(x), \\
& y_2\vert_{\ppp\OOO\times (0,T)} = 0,
\end{array}\right.
                              \eqno{(4.18)}
$$
where $\www{S_1}$ satisfies (4.14).  Similarly to (4.16) and (4.17), 
we can verify
$$
E(y_2)(0) =: E(0) \le \Vert a\Vert^2_{H^3(\OOO)}, \quad
\Vert \www{S_1}\Vert^2_{L^2(0,T;L^2(\OOO))} \le CE(0).
$$
Applying (3.20) in Lemma 3.3, we have
$$
\Vert \ppp_t^2\ppp_{\nu}y\Vert^2_{L^2(0,T;L^2(\ppp\OOO))}
\le C(E(0) + \Vert \www{S_1}\Vert^2_{L^2(0,T;L^2(\OOO))})  
                                                  \eqno{(4.19)}
$$
$$
\le CE(0) \le C\Vert a\Vert^2_{H^3(\OOO)}.
$$
By $\ppp_{\nu}\ppp_ty = 0$ on $\ppp\OOO$, we see
$$
\ppp_t\ppp_{\nu}y(x,t) = \int^t_0 \ppp_t^2\ppp_{\nu}y(x,\eta)d\eta
$$
and
$$
\ppp_{\nu}y(x,t) = \int^t_0 \ppp_t\ppp_{\nu}y(x,\eta)d\eta
+ \ppp_{\nu}a.
$$
Therefore
$$
\Vert \ppp_t\ppp_{\nu}y\Vert_{L^2(0,T;L^2(\ppp\OOO))}
\le C\Vert \ppp_t^2\ppp_{\nu}y\Vert_{L^2(0,T;L^2(\ppp\OOO))}
$$
and
\begin{align*}
&\Vert \ppp_{\nu}y\Vert_{L^2(0,T;L^2(\ppp\OOO))}
\le C\Vert \ppp_t^2\ppp_{\nu}y\Vert_{L^2(0,T;L^2(\ppp\OOO))}
+ C\Vert \ppp_{\nu}a\Vert_{L^2(\ppp\OOO)}\\
\le & C\Vert \ppp_t^2\ppp_{\nu}y\Vert_{L^2(0,T;L^2(\ppp\OOO))}
+ C\Vert a\Vert_{H^3(\OOO)}.
\end{align*}
With (4.19), we complete the proof of the first inequality.
Thus the proof of Theorem 1.2 is completed.
\section{Proof of Theorem 1.3}

Once that the Carleman estimate Theorem 1.1 is established, we can 
prove Theorem 1.3 by an argument similar to Imanuvilov and Yamamoto 
\cite{IY3}.  See also Bellassoued and Yamamoto \cite{BelY3}.
\\
{\bf First Step.}\\
We set $F(x,t) = R(x,t)f(x)$ and
$$
v(x,t) = \chi Au(x,t) + \chi(t)\int^t_0 B(t,\eta)u(x,\eta) d\eta,
\quad (x,t) \in \Omega \times (0,T).            \eqno{(5.1)}
$$
Similarly to (2.8) in $\OOO \times (0,T)$, using $f = \ppp_{\nu}f = 0$
on $\ppp\OOO$, by (2.4) we can verify
$$
v(x,t) = \chi(t)\ppp_t^2u(x,t) - \chi(t)F(x,t), \quad 
x\in \OOO, \thinspace 0 < t < T
$$
and
$$\left\{ \begin{array}{rl}
& \ppp_t^2v(x,t) - Av = \chi AF + \chi J + S, \quad x \in \OOO, \thinspace
0 < t < T,\\
& v(x,0) = \ppp_tv(x,0) = 0, \quad x \in \OOO,\\
& \ppp_t^jv(\cdot,T) = 0 \quad \mbox{in $\OOO$, $j=0,1$},\\
& v\vert_{\ppp\OOO\times (0,T)} = 0.
\end{array}\right.           \eqno{(5.2)}
$$
We make the even extension of $v$ to $(-T,0)$:
$$
v(x,t) =
\left\{ \begin{array}{rl}
v(x,t), \quad & 0 < t < T, \\
v(x,-t), \quad & -T < t < 0.
\end{array}\right.
$$
Accordingly we make the even extensions of $\chi AF + \chi J + S$.  Then,
by $v(x,0) = \ppp_tv(x,0) = 0$ for $x \in \OOO$ we can prove that 
$$
v \in C^2([-T,T];L^2(\Omega)) \cap C^1([-T,T];H^1_0(\Omega))
\cap C([-T,T]; H^2(\OOO)\cap H^1_0(\OOO)) 
$$
$$
\cap H^3(-T,T;L^2(\OOO))      
                                           \eqno{(5.3)}
$$
and
$$
\chi AF + \chi J + S \in H^1(-T,T;L^2(\OOO)).    \eqno{(5.4)}
$$
We recall that $Q = \OOO \times (-T,T)$.
Hence, setting $v_1 = \ppp_tv$, we have
$$\left\{ \begin{array}{rl}
& \ppp_t^2v_1 - Av_1 = \chi A\ppp_tF + \chi'AF + \chi(\ppp_tJ)
+ \chi'J + \ppp_tS \quad \mbox{in $Q$}, \\
& \ppp_t^jv_1(\cdot,\pm T) = 0 \quad \mbox{in $\OOO$, $j=0,1$},\\
& v_1\vert_{\ppp\OOO} = 0,\\
& v_1(x,0) = 0, \quad \ppp_tv_1(x,0) = \chi AF(x,0).
\end{array}\right.           \eqno{(5.5)}
$$
Here we used that $\ppp_tv_1(x,0) = \chi AF(x,0)$ by $\ppp_tv_1
= \ppp_t^2v$ and (2.7), (2.9).
 
We set 
$$
z = (\ppp_tv)e^{s\va} = v_1e^{s\va}.
$$
We write (5.5) in terms of $z$.  First we have
$$
\ppp_tz = (\ppp_tv_1)e^{s\va} + s(\ppp_t\va)v_1e^{s\va}
$$
and
$$
\ppp_t^2 z = (\ppp_t^2v_1)e^{s\va} + 2s(\ppp_t\va)(\ppp_tv_1)e^{s\va}
+ s(\ppp_t^2\va)v_1e^{s\va} + s^2(\ppp_t\va)^2v_1e^{s\va}.
$$
Moreover 
$$
\ppp_iz = (\ppp_iv_1)e^{s\va} + s(\ppp_i\va)v_1e^{s\va},
$$
and so
\begin{align*}
& \ppp_j\ppp_i z = (\ppp_j\ppp_iv_1)e^{s\va} 
+ s\{ (\ppp_i\va)(\ppp_jv_1) + (\ppp_j\va)(\ppp_iv_1)\}e^{s\va}\\
+& \{s(\ppp_i\ppp_j\va) + s^2(\ppp_i\va)\ppp_j\va\}v_1e^{s\va}.
\end{align*}
Hence 
\begin{align*}
&Az = \sumij a_{ij}\ppp_i\ppp_jz + \sumij (\ppp_ia_{ij})\ppp_jz\\
=& \sumij a_{ij}(\ppp_i\ppp_jv_1)e^{s\va}
+ \sumij sa_{ij}((\ppp_iv_1)\ppp_j\va + (\ppp_jv_1)\ppp_i\va)e^{s\va}\\
+ & \sumij sa_{ij}(\ppp_i\ppp_j\va)v_1 e^{s\va}
+ s^2\sumij a_{ij}(\ppp_i\va)(\ppp_j\va)v_1e^{s\va}\\
+& \sumij (\ppp_ia_{ij})(\ppp_jv_1)e^{s\va} 
+ \sumij (\ppp_ia_{ij})s(\ppp_j\va)v_1e^{s\va}.
\end{align*}
Using $a_{ij} = a_{ji}$ in the second term on the right-hand side, we
obtain
\begin{align*}
& Az = e^{s\va}Av_1 + 2s\sumij a_{ij}(\ppp_i\va)(\ppp_jv_1)e^{s\va}\\
+& sv_1e^{s\va}A\va + s^2v_1e^{s\va}\sumij a_{ij}(\ppp_i\va)\ppp_j\va.
\end{align*}
Thus (5.5) yields
$$
\ppp_t^2 z - Az                    \eqno{(5.6)}
$$
\begin{align*}
= &\chi e^{s\va}A\ppp_tF + \chi'e^{s\va}AF
+ \chi(\ppp_tJ)e^{s\va} + \chi'Je^{s\va} + (\ppp_tS)e^{s\va}\\
+ & 2s((\ppp_t\va)\ppp_tv_1 - \sumij a_{ij}(\ppp_iv_1)\ppp_j\va)e^{s\va}
+ s(\ppp_t^2\va - A\va)v_1e^{s\va}\\
+& s^2((\ppp_t\va)^2 - \sumij a_{ij}(\ppp_i\va)\ppp_j\va)v_1e^{s\va}
\quad \mbox{in $Q$},
\end{align*}
$$\left\{ \begin{array}{rl}
& \ppp_tz(x,0) = \chi A(R(x,0)f(x))e^{s\va(x,0)}, \\
& z(x,0) = 0, \qquad x \in \OOO,
\end{array}\right.                          \eqno{(5.7)}
$$
and
$$
z\vert_{\ppp\OOO} = 0.                   \eqno{(5.8)}
$$

We rewrite Lemma 2.4 in terms of $z=v_1e^{s\va}$.  First $s^3\vert z\vert^2
= s^3\vert v_1\vert^2\weight$ and then
$\ppp_tz = (\ppp_tv_1)e^{s\va} + s(\ppp_t\va)z$, so that
$$
s\vert \ppp_t z\vert^2 \le Cs^3\vert z\vert^2 + Cs\vert \ppp_tv_1\vert^2
\weight,
$$
and we have similar estimates for $s\vert \nabla z\vert^2$.
Hence
$$
\int_Q (s\vert \nabla_{x,t}z\vert^2 + s^3\vert z\vert^2) dxdt
\le C\int_Q (s\vert \nabla_{x,t}v_1\vert^2 + s^3\vert v_1\vert^2)\weight
dxdt,
$$
and so
\begin{align*}
& \int_Q (s\vert \nabla_{x,t}z\vert^2 + s^3\vert z\vert^2) dxdt
\le C\int_Q (\vert AF\vert^2 + \vert A\ppp_tF\vert^2)\weight dxdt\\
+& C\Vert u\Vert^2_{H^2(-T,T;H^2(\OOO))}s^2e^{2s\delta}
+ Ce^{Cs}\Vert \ppp_{\nu}u\Vert^2_{H^3(-T,T;L^2(\Gamma))}.
\end{align*}
for $s > s_0$.  Here we include $\Phi^2$ into $C>0$.
Henceforth we set 
$$
D^2 = \Vert \ppp_{\nu}u\Vert^2_{H^3(-T,T;L^2(\Gamma))}.
$$

We estimate $\Vert u\Vert^2_{H^2(-T,T;H^2(\OOO))}$ by Lemma 3.2, so that 
we obtain
$$
\int_Q (s\vert \nabla_{x,t}z\vert^2 + s^3\vert z\vert^2) dxdt
\le C\int_Q (\vert AF\vert^2 + \vert A\ppp_tF\vert^2)\weight dxdt
                                   \eqno{(5.9)}
$$
$$
+  C\Vert f\Vert^2_{H^2(\OOO)}s^2e^{2s\delta} + Ce^{Cs}D^2
$$
for $s > s_0$.  
\\
{\bf Second Step.}\\
We will carry out the energy estimate.  We multiply (5.6) 
with $\ppp_tz$ and integrate by parts over $\OOO \times (-T,0)$. Then
\begin{align*}
& \int^0_{-T}\int_{\OOO} (\ppp_t^2z - Az)\ppp_tz dxdt\\
= & \frac{1}{2}\int^0_{-T} \ppp_t\left(\int_{\OOO} \vert \ppp_tz\vert^2 dx
\right) dt
- \int^0_{-T}\int_{\OOO} (Az)\ppp_tz dxdt.
\end{align*}
Here by $z\vert_{\ppp\OOO} = 0$ and $a_{ij}=a_{ji}$, integrating by parts,
we see
\begin{align*}
&-\int_{\OOO} (Az)\ppp_tz dx = \sumij \int_{\OOO} 
a_{ij}(\ppp_i\ppp_tz)\ppp_jz dx\\
=& \int_{\OOO} \biggl\{ \sum_{i<j} a_{ij}((\ppp_iz)(\ppp_j\ppp_tz)
+ (\ppp_jz)(\ppp_i\ppp_tz))
+ \sum_{i=1}^n a_{ii}(\ppp_iz)(\ppp_i\ppp_tz) \biggr\} dx\\
=& \int_{\OOO} \left\{ \sum_{i>j} a_{ij}\ppp_t((\ppp_iz)\ppp_jz)
+ \sum_{i=1}^n a_{ii}\frac{1}{2} \ppp_t(\vert \ppp_iz\vert^2) \right\}
dx\\
= & \frac{1}{2} \int_{\OOO} \sumij a_{ij}\ppp_t((\ppp_iz)\ppp_jz)dx.
\end{align*}
Using $z(\cdot,0) = z(\cdot, -T) = 0$ in $\OOO$ and noting that 
$a_{ij}$ is independent of $t$, we obtain
$$
\int^0_{-T}\int_{\OOO} (\ppp_t^2z - Az)\ppp_tz dxdt
= \frac{1}{2}\int_{\OOO} \vert \ppp_tz(x,0)\vert^2 dx.
                                                      \eqno{(5.10)}
$$

Next, by the Cauchy-Schwarz inequality, we have
$$
\left\vert \int^0_{-T}\int_{\OOO} \mbox{[the right-hand side] of (5.6)]}
\times \ppp_tz dxdt \right\vert                   \eqno{(5.11)}
$$
\begin{align*}
= &\biggl\vert \int^0_{-T} \int_{\OOO}
(\chi A\ppp_tF + \chi'AF) \weight \ppp_tz dxdt
+ \int^0_{-T} \int_{\OOO} (\chi(\ppp_tJ) + \chi'J + \ppp_tS)\weight 
\ppp_tz dxdt\\
+ & 2s\int^0_{-T} \int_{\OOO}
((\ppp_t\va)\ppp_tv_1 - \sumij a_{ij}(\ppp_iv_1)\ppp_j\va)
\weight \ppp_tz dxdt\\
+& s \int^0_{-T} \int_{\OOO} (\ppp_t^2\va-A\va)v_1\weight 
\ppp_tz dxdt\\
+& s^2\int^0_{-T} \int_{\OOO}((\ppp_t\va)^2 - \sumij a_{ij}(\ppp_i\va)
(\ppp_j\va))v_1\weight \ppp_tz dxdt\biggr\vert\\
\le& C\int_Q \vert A\ppp_tF\vert^2 + \vert AF\vert^2) \weight dxdt
+ C\int_Q \vert \ppp_tz\vert^2 dxdt\\
+& C\int_Q (\vert \ppp_tJ\vert^2 + \vert J\vert^2 + \vert \ppp_tS\vert^2)
\weight dxdt\\
+& C\int_Q s\vert \nabla_{x,t}v_1\vert e^{s\va} \vert \ppp_tz\vert dxdt
+ Cs^2\int_Q \vert v_1\vert e^{s\va}\vert \ppp_tz\vert dxdt.
\end{align*}
Here we extended the integral domain $\OOO \times (-T,0)$ to $\OOO \times 
(-T,T) =: Q$.  Moreover we note
$$
\int_Q s\vert \nabla_{x,t}v_1\vert e^{s\va} \vert \ppp_tz\vert dxdt
\le 2\int_Q  s(\vert \nabla_{x,t}v_1\vert^2 \weight 
+ \vert \ppp_tz\vert^2) dxdt
$$
and
\begin{align*}
&\int_Q s^2\vert v_1\vert e^{s\va}\vert \ppp_tz\vert dxdt
= \int_Q s^{\frac{3}{2}}\vert v_1\vert e^{s\va}
s^{\frac{1}{2}}\vert \ppp_tz\vert dxdt\\
\le& 2\int_Q (s^3\vert v_1\vert^2 e^{2s\va} + s\vert \ppp_tz\vert^2)
dxdt,
\end{align*}
By (2.22) and (1.13), we have
\begin{align*}
& \int_Q (\vert \ppp_tJ\vert^2 + \vert J\vert^2) \weight dxdt
\le C\int_Q \sumalpha \sum_{k=0}^2 \vert \dalpha(\chi \ppp_t^ku)\vert^2
\weight dxdt\\
\le& C\int_Q (\vert A\ppp_tF\vert^2 + \vert AF\vert^2) \weight dxdt
+ C\Vert f\Vert^2_{H^2(\OOO)}s^2e^{2s\delta} + Ce^{Cs}D^2,
\end{align*}
and (2.23) yields
\begin{align*}
&\int_Q \vert \ppp_tS\vert^2 \weight dxdt 
\le C\Vert u\Vert_{H^2(-T,T;H^2(\OOO))}^2s^2e^{2s\delta}\\
\le &C\Vert f\Vert^2_{H^2(\OOO)}s^2e^{2s\delta} + Ce^{Cs}D^2.
\end{align*}
Consequently substituting these inequalities and applying
Lemma 2.4 and (5.9) to estimate
$$
\int_Q (\vert \ppp_tz\vert^2 + s\vert \nabla_{x,t}v_1\vert^2\weight
+ s\vert \ppp_tz\vert^2 + s^3\vert v_1\vert^2\weight) dxdt,
$$
from (5.11) we reach
$$
\left\vert \int^0_{-T} \int_{\OOO}
\mbox{[the right-hand side of (5.6)]}\thinspace \times \ppp_tz dxdt
\right\vert                                                \eqno{(5.12)}
$$
$$
\le C\int_Q (\vert A\ppp_tF\vert^2 + \vert AF\vert^2)\weight dxdt
+ C\Vert f\Vert^2_{H^2(\OOO)}s^2e^{2s\delta} + Ce^{Cs}D^2
$$
for $s>s_0$.

By (5.7), (5.10) and (5.12), using $\chi(0) = 1$, we see
\begin{align*}
& \int_{\OOO} \vert \ppp_tz(x,0)\vert^2 dx
= \int_{\OOO} \vert A(R(x,0)f(x))\vert^2 \weight dxdt\\
\le & C\int_Q (\vert A\ppp_tF\vert^2 + \vert AF\vert^2)\weight dxdt
+ C\Vert f\Vert^2_{H^2(\OOO)}s^2e^{2s\delta} + Ce^{Cs}D^2
\end{align*}
for $s > s_0$.
\\
{\bf Third Step.}\\
We complete the proof by an elliptic Carleman estimate.
Since
\begin{align*}
& A((\ppp_t^kR)f) = (\ppp_t^kR)(Af) + \sumij (\ppp_i\ppp_t^kR)
a_{ij}(\ppp_jf)\\
+& \sumij \ppp_i(a_{ij}(\ppp_j\ppp_tR)f), \quad k=0,1,
\end{align*}
we estimate
$$
\vert A(R(x,0)f)\vert \ge \vert R(x,0)Af\vert 
- C(\vert \nabla f\vert + \vert f\vert)
$$
and
$$
\vert A((\ppp_t^kR)f)\vert^2\le C(\vert Af\vert^2 + \vert \nabla f\vert^2
+ \vert f\vert^2),
$$
and so
\begin{align*}
& \int_{\OOO} \vert R(0)Af\vert^2 e^{2s\va(x,0)} dx
- C\int_{\OOO} (\vert \nabla f\vert^2 + \vert f\vert^2)
e^{2s\va(x,0)} dx\\
\le &C\int_{\OOO}\left( \int^T_{-T} \vert Af\vert^2 \weight dt\right) dx
+ C\int_{\OOO}\int^T_{-T} (\vert \nabla f\vert^2 + \vert f\vert^2)
\weight dxdt \\
+& C\Vert f\Vert^2_{H^2(\OOO)}s^2e^{2s\delta} + Ce^{Cs}D^2.
\end{align*}
By $\va(x,t) \le \va(x,0)$, we replace the second term on the right-hand side
by the second term on the left-hand side, and we apply $\vert R(\cdot,0)\vert
> 0$ on $\overline{\OOO}$ by (1.17).
Therefore 
$$
\int_{\OOO} \vert Af\vert^2 e^{2s\va(x,0)} dx
\le C\int_{\OOO} \vert Af\vert^2 e^{2s\va(x,0)} 
\left( \int^T_{-T} r^{2s(\va(x,t) - \va(x,0))} dt \right)dx 
                               \eqno{(5.13)}
$$
$$
+ C\int_{\OOO} (\vert \nabla f\vert^2 + \vert f\vert^2)
e^{2s\va(x,0)} dx 
+ C\Vert f\Vert^2_{H^2(\OOO)}s^2e^{2s\delta} + Ce^{Cs}D^2.
$$
Since 
\begin{align*}
& 2s(\va(x,t)-\va(x,0)) 
= 2s(e^{\gamma\vert x-x_0\vert^2 - \gamma\beta t^2}
- e^{\gamma\vert x-x_0\vert^2})\\
=& 2se^{\gamma\vert x-x_0\vert^2}(e^{-\gamma\beta t^2}-1)
\le 2s(e^{-\gamma\beta t^2} - 1)
\end{align*}
and $e^{-\gamma\beta t^2} - 1 < 0$ for $t \ne 0$, the Lebesgue theorem 
yields
$$
\int^T_{-T} e^{2s(\va(x,t) - \va(x,0))} dt
\le \int^T_{-T} e^{2s(e^{-\gamma\beta t^2} - 1)} dt = o(1)
$$
as $s \to \infty$.

Therefore by choosing $s>0$ sufficiently large, the estimate 
(5.13) implies
\begin{align*}
& \frac{1}{2}\int_{\OOO} \vert Af\vert^2 e^{2s\va(x,0)} dx
\le (1-o(1))\int_{\OOO} \vert Af\vert^2 e^{2s\va(x,0)} dx\\
\le& C\int_{\OOO} (\vert \nabla f\vert^2 + \vert f\vert^2)
e^{2s\va(x,0)} dx
+ C\Vert f\Vert^2_{H^2(\OOO)}s^2e^{2s\delta} + Ce^{Cs}D^2.
\end{align*}
for $s > s_0$.
By $\nabla \va(x,0) = 2\gamma(x-x_0)\va \ne 0$ for $x \in \ooo{\OOO}$,
we apply the Carleman estimate for the elliptic operator $A$ of the second
order which is similar to Lemma 2.2 (here we fix $\gamma$), we have
\begin{align*}
& \int_{\OOO} \left( \frac{1}{s}\sum_{\vert \alpha\vert = 2}
\vert \dalpha f\vert^2 + s\vert \nabla f\vert^2 + s^3\vert f\vert^2\right)
e^{2s\va(x,0)} dx\\
\le &C\int_{\OOO} (\vert \nabla f\vert^2 + \vert f\vert^2)
e^{2s\va(x,0)} dx 
+ C\Vert f\Vert^2_{H^2(\OOO)}s^2e^{2s\delta} + Ce^{Cs}D^2
\end{align*}
for $s > s_0$.  Again choosing $s>0$ sufficienly large, we can 
absorb the first term on the right-hand side into the left-hand side,  
multiplying with $s$ and replacing $se^{Cs}$ by $e^{Cs}$,
we have
$$
\int_{\OOO} \sumalpha \vert \dalpha f\vert^2 e^{2s\va(x,0)}dx
\le C\Vert f\Vert^2_{H^2(\OOO)}s^3e^{2s\delta}
+ Ce^{Cs}D^2                                 \eqno{(5.14)}
$$
for $s > s_1$.  

Hence
$$
e^{2s(\delta+2\ep_0)}\int_{\OOO} \sum_{\vert \alpha\vert \le 2}
\vert \dalpha f\vert^2 dx
\le C\Vert f\Vert^2_{H^2(\OOO)}s^3e^{2s\delta}
+ Ce^{Cs}D^2,
$$
that is,
$$
(1 - Cs^3e^{-4\ep_0s})\Vert f\Vert^2_{H^2(\OOO)}
\le Ce^{Cs}D^2
$$
for all $s > s_0$.
We choose $s>0$ sufficiently large so that 
$1 - Cs^3e^{-4\ep_0s} > 0$.
Then $\Vert f\Vert^2_{H^2(\OOO)} \le C_1e^{Cs}D^2$.
Noting that 
$$
D^2 = \Vert \ppp_{\nu}u\Vert^2_{H^3(-T,T;L^2(\Gamma))}
= 2\Vert \ppp_{\nu}u\Vert^2_{H^3(0,T;L^2(\Gamma))},
$$
we complete the proof of the second inequality of (1.18).

Finally we have to prove the first inequality in (1.18).
We set $w = \ppp_t^3u$.  Similarly to (2.8), we can obtain
$$\left\{ \begin{array}{rl}
& \ppp_t^2 w(x,t) = Aw + S_2(x,t) + (\ppp_t^3R)f, \quad x\in \OOO,
\thinspace 0 < t < T, \\
& w(x,0) = (\ppp_tR)(x,0)f(x), \\
& \ppp_tw(x,0) = A(R(x,0)f)(x) + (\ppp_t^2R)(x,0)f(x), \quad 
x\in \OOO, \\
& w\vert_{\ppp\OOO\times (0,T)} = 0,
\end{array}\right.
                              \eqno{(5.15)}
$$
where $\www{B}(t) = B(t,t)$ and 
\begin{align*}
& S_2(x,t) = \{(\ppp_t^2\www{B})(t) + \ppp_t((\ppp_tB)(t,t))
+ (\ppp_t^2B)(t,t)\} u(x,t)\\
+& \{ 2(\ppp_t\www{B})(t) + (\ppp_tB)(t,t)\}\ppp_tu(x,t)\\
+& \www{B}(t)\ppp_t^2u(x,t) + \int^t_0 \ppp_t^3B(t,\eta)u(x,\eta)d\eta.
\end{align*}
Therefore by (1.4) we see 
$$
\vert S_2(x,t)\vert \le C\left( \sumalpha \sum_{k=0}^2
\vert \dalpha \ppp_t^ku(x,t)\vert
+ \int^t_0 \sumalpha \vert \dalpha u(x,\eta)\vert d\eta\right)
                                                  \eqno{(5.16)}
$$
for $x\in \OOO$ and $0 < t < T$.
We set 
$$
E(t) = E(w)(t) = \int_{\OOO} (\vert \ppp_tw(x,t)\vert^2
+ \sumij a_{ij}\ppp_jw(x,t)\ppp_iw(x,t)) dx.
$$
Then we readily verify that $E(0) \le C\Vert f\Vert^2_{H^2(\OOO)}$.
Applying (3.20) to (5.15) and noting (5.16), we obtain
\begin{align*}
& \Vert \ppp_{\nu}w\Vert^2_{L^2(0,T;L^2(\ppp\OOO))}
= \Vert \ppp_t^3\ppp_{\nu}u\Vert^2_{L^2(0,T;L^2(\ppp\OOO))}\\
\le& C(\Vert f\Vert^2_{H^2(\OOO)} + \Vert S_2\Vert^2
_{L^2(0,T;L^2(\OOO))})\\
\le& C\left(\Vert f\Vert^2_{H^2(\OOO)} + \Vert u\Vert^2_{H^2(0,T;H^2(\OOO))}
+ \int^T_0 \left( \int^t_0 \Vert u(\cdot,\eta)\Vert^2_{H^2(\OOO)}d\eta
\right) dt \right)\\
\le& C(\Vert f\Vert^2_{H^2(\OOO)} + \Vert u\Vert^2
_{H^2(0,T;H^2(\OOO))}).
\end{align*}
The second term on the right-hand side is estimated by Lemma 3.2, so that
$$
\Vert \ppp_t^3\ppp_{\nu}\Vert^2_{L^2(0,T;L^2(\ppp\OOO))}
= C\Vert f\Vert^2_{H^2(\OOO)}.
$$
Finally, since
$$
\ppp_t^2\ppp_{\nu}u(x,0) = \ppp_{\nu}(R(x,0)f) = 0,
\quad \ppp_t\ppp_{\nu}u(x,0) = 0
$$
by $f \in H^2_0(\OOO)$ and $u(x,0) = \ppp_tu(x,0) = 0$, we have
$$
\ppp_t^k\ppp_{\nu}u(x,t) = \int^t_0 \ppp_t^{k+1}\ppp_{\nu}(x,\xi) d\xi,
\quad k=0,1,2.
$$
Consequently $\Vert \ppp_{\nu}u\Vert_{H^3(0,T;L^2(\ppp\OOO))}
\le C\Vert f\Vert_{H^2(\OOO)}$.  Thus the proof of the first inequality, and
so Theorem 1.3 is completed.
\\
\vspace{0.2cm}

{\bf Acknowledgements.}  This work was completed when the third author
was a guest professor at Sapienza Universit\`a di Roma in May - June 2016.
The author thanks the university for that opportunity yielding the current 
joint work. 
\\
\vspace{0.2cm}

{\bf Appendix. Proof of Lemma 2.2 for general $p$.}

It is known that
$$
\int_Q \left(s^{-1}\va^{-1} \sumalalpha \vert \ppp_x^{\alpha}y\vert^2
+ s\gamma^2\va\vert \nabla y\vert^2 + s^3\gamma^4\va^3\vert y\vert^2
\right) \weight dxdt                          \eqno{(1)}
$$
$$
\le C\int_Q \vert Ay\vert^2 \weight dxdt + Ce^{Cs}\Vert \ppp_{\nu}y\Vert^2
_{L^2(\Sigma)}
$$
for all $s > s_0$ and $y \in L^2(-T,T;H^2(\Omega) \cap H^1_0(\Omega))$.
See e.g., Bellassoued and Yamamoto \cite{BelY3} or we can prove (1) 
similarly to Yamamoto \cite{Y} where the parabolic Carleman estimate
is proved.  However here we omit the proof of (1).
Now let $p \ne -1$, $\in \R$.  We set 
$$
\theta = \frac{p+1}{2} \ne 0             \eqno{(2)}
$$
and
$$
z = y\va^{\theta}, \quad Ay = F.          \eqno{(3)}
$$
Then we directly verify 
$$
\ppp_j\va = \gamma\va(\ppp_j\psi), \quad
\ppp_i\ppp_j\va = \gamma(\ppp_i\ppp_j\psi)\va 
+ \gamma^2(\ppp_i\psi)(\ppp_j\psi)\va,
$$
and
$$
\ppp_jy = -\gamma\theta\va^{-\theta}(\ppp_j\psi)z + \va^{-\theta}
\ppp_jz                      \eqno{(4)}
$$
and
$$
\ppp_i\ppp_jy = \gamma^2\theta^2\va^{-\theta}(\ppp_i\psi)(\ppp_j\psi)z
- \gamma\theta\va^{-\theta}(\ppp_i\ppp_j\psi)z
$$
$$
- \gamma\theta\va^{-\theta}((\ppp_j\psi)\ppp_iz + (\ppp_i\psi)\ppp_jz)
+ \va^{-\theta}\ppp_i\ppp_jz, \quad 1\le j \le n.   
                                             \eqno{(5)}
$$
Therefore we have
\begin{align*}
& Ay = \va^{-\theta}Az 
- 2\theta \gamma \va^{-\theta}\sum_{i,j=1}^n a_{ij}(\ppp_i\psi)(\ppp_jz)\\
+& \va^{-\theta}\left\{ \sum_{i,j=1}^n a_{ij}\theta^2\gamma^2
(\ppp_i\psi)\ppp_j\psi
- \theta a_{ij}\gamma(\ppp_i\ppp_j\psi)
- \theta\gamma (\ppp_ia_{ij})\ppp_j\psi\right\}z \quad \mbox{in $Q$},
\end{align*}
that is,
\begin{align*}
& Az = \va^{\theta}F
+ 2\theta \gamma \sum_{i,j=1}^n a_{ij}(\ppp_i\psi)(\ppp_jz)\\
-& \left\{ \sum_{i,j=1}^n a_{ij}\theta^2\gamma^2
(\ppp_i\psi)\ppp_j\psi
- \theta a_{ij}\gamma(\ppp_i\ppp_j\psi)
- \theta\gamma (\ppp_ia_{ij})\ppp_j\psi\right\}z \quad \mbox{in $Q$}.
\end{align*}
By $a_{ij} \in C^1(\overline{\Omega})$, we can write
$$
\vert Az\vert^2 \le C(\va^{2\theta}\vert F\vert^2
+ \gamma^2\vert \nabla z\vert^2 + \gamma^4\vert z\vert^2) \quad 
\mbox{in $Q$}.                               \eqno{(6)}
$$
By (3), applying (1) to (6), we obtain
\begin{align*}
& \int_Q \left(s^{-1}\va^{-1} \sumalalpha \vert \ppp_x^{\alpha}y\vert^2
+ s\gamma^2\va\vert \nabla z\vert^2 + s^3\gamma^4\va^3\vert z\vert^2
\right)\weight dxdt\\
\le& C\int_Q \va^{p+1}\vert F\vert^2 \weight dxdt 
+ \int_Q (\gamma^2\vert \nabla z\vert^2 + \gamma^4\vert z\vert^2)
\weight dxdt \\
+& Ce^{Cs}\Vert \ppp_{\nu}z\Vert^2_{L^2(\Sigma)}.
\end{align*}
Choosing $s>0$ sufficiently large, we can absorb the second term on the 
right-hand side into the left-hand side and using
$$
\ppp_{\nu}z = \nabla (y\va^{\theta})\cdot \nu = \va^{\theta}\ppp_{\nu}y
\quad \mbox{on $\ppp\Omega$}
$$
by $y\vert_{\ppp\Omega} = 0$, we have
$$
\int_Q \left(s^{-1}\va^{-1} \sumalalpha \vert \ppp_x^{\alpha}y\vert^2
+ s\gamma^2\va\vert \nabla z\vert^2 + s^3\gamma^4\va^3\vert z\vert^2
\right) \weight dxdt
$$
$$
\le C\int_Q \va^{p+1}\vert F\vert^2 \weight dxdt 
+ Ce^{Cs}\Vert \va^{p+1}\Vert_{L^{\infty}(\Sigma)}
\Vert \ppp_{\nu}y\Vert^2_{L^2(\Sigma)}.               \eqno{(7)}
$$
Similarly to (4) and (5), we see
$$
\ppp_jz = \gamma\theta\va^{\theta}(\ppp_j\psi)y + \va^{\theta}
\ppp_jy
$$
and
\begin{align*}
& \ppp_i\ppp_jz = \gamma^2\theta^2\va^{\theta}(\ppp_i\psi)(\ppp_j\psi)y
+ \gamma\theta\va^{\theta}(\ppp_i\ppp_j\psi)y\\
+ &\gamma\theta\va^{\theta}((\ppp_j\psi)\ppp_iy + (\ppp_i\psi)\ppp_jy)
+ \va^{\theta}\ppp_i\ppp_jy, \quad 1\le j \le n.
\end{align*}
Therefore, since 
$$
s^3\gamma^4\va^3\va^{p+1}\vert y\vert^2 = s^3\gamma^4\va^3\vert z\vert^2,   
                                               \eqno{(8)}
$$
we obtain
$$
s\gamma^2\va^{p+2}\vert \nabla y\vert^2
\le C(s\gamma^2\va\vert \nabla z\vert^2 + s\gamma^4\va^{p+2}\vert y\vert^2)
$$
$$
\le C(s\gamma^2\va\vert \nabla z\vert^2 + s^3\gamma^4\va^3\vert z\vert^2).
                                                 \eqno{(9)}
$$
Here we used (8) for the final term.
Moreover by (8) and (9), we have
\begin{align*}
&s^{-1}\va^{-1}\va^{p+1}\vert \ppp_i\ppp_jy\vert^2
\le Cs^{-1}\va^{-1}(\gamma^4\va^{p+1}\vert y\vert^2 
+ \gamma^2\va^{p+1}\vert \nabla y\vert^2 + \vert \ppp_i\ppp_jz\vert^2)\\
\le& C(s^{-1}\va^{-1}\vert \ppp_i\ppp_jz\vert^2 
+ s^{-1}\gamma^2\va^p \vert \nabla y\vert^2 + s^{-1}\gamma^4\va^p
\vert y\vert^2)
\end{align*}
$$
\le C(s^{-1}\va^{-1}\vert \ppp_i\ppp_jz\vert^2 
+ s\gamma^2\va \vert \nabla z\vert^2 + s^3\gamma^4\va^3\vert z\vert^2).
                                                \eqno{(10)}
$$
Substituting (8) - (10) in (7), we obtain
$$
\int_Q \left(s^{-1}\va^p \sumalalpha \vert \ppp_x^{\alpha}y\vert^2
+ s\gamma^2\va^{p+2}\vert \nabla y\vert^2 + s^3\gamma^4\va^{p+4}\vert y\vert^2
\right) \weight dxdt                          \eqno{(11)}
$$
$$
\le C\int_Q \va^{p+1}\vert F\vert^2 \weight dxdt 
+ Ce^{Cs}e^{C\gamma\Vert \psi\Vert_{L^{\infty}(\Sigma)}}
\Vert \ppp_{\nu}y\Vert^2
_{L^2(\Sigma)}
$$
for all $s > s_0$.  Here in choosing $s_0(\gamma)$, we further 
assume that $s_0(\gamma) = s_0 > \gamma$ to have 
$$
s^{p+1}e^{Cs}e^{C\gamma\Vert \psi\Vert_{L^{\infty}(\Sigma)}}
\le e^{C_1s}
$$
with sufficiently large constant $C_1 > 0$.  Multiplying (11) with 
$s^{p+1}$, we reach the conclusion and thus the proof of Lemma 2.2
with $p \in \R$ is completed.

\end{document}